\documentclass[12pt,oneside,reqno]{amsart}

\usepackage[UKenglish]{babel}
\usepackage[a4paper,top=2cm,bottom=2cm,left=3cm,right=3cm]{geometry}

\usepackage{mathtools}
\usepackage{amssymb}
\usepackage{bbm}
\usepackage{dsfont}
\usepackage{scalerel}
\usepackage{stackengine}
\usepackage{cancel}

\usepackage{graphicx}
\usepackage{tikz}
\usepackage{tikz-cd}

\usepackage{enumerate}
\usepackage{setspace}
\usepackage{tabularx}
\usepackage{xcolor}
\usepackage[normalem]{ulem}
\usepackage{contour}
\usepackage{dirtytalk}

\usepackage[breaklinks]{hyperref}
\newcounter{eqH}
\makeatletter
\AtBeginDocument{%
	\let\HRO@refstepcounter\refstepcounter
	\renewcommand{\refstepcounter}[1]{%
		\def\HRO@tmpa{#1}\def\HRO@tmpb{equation}%
		\ifx\HRO@tmpa\HRO@tmpb\stepcounter{eqH}\fi
		\HRO@refstepcounter{#1}}%
}
\makeatother

\DeclareRobustCommand{\myuline}[1]{%
	\uline{\phantom{#1}}%
	\llap{\contour{white}{#1}}%
}

\theoremstyle{plain}
\newtheorem{theorem}{Theorem}[section]
\newtheorem{lemma}[theorem]{Lemma}
\newtheorem{corollary}[theorem]{Corollary}

\theoremstyle{definition}

\theoremstyle{remark}
\newtheorem{remark}[theorem]{Remark}

\newcommand{\f}{\frac}
\renewcommand{\i}{\infty}
\renewcommand{\t}{\text}

\renewcommand{\l}{\left}
\renewcommand{\r}{\right}

\newcommand{\littletaller}{\mathchoice{\vphantom{\big|}}{}{}{}}
\newcommand\restr[2]{{
		\left.\kern-\nulldelimiterspace 
		#1 
		\littletaller 
		\right|_{#2} 
}}

\newcommand{\om}{{\omega}}

\newcommand{\la}{\lambda}

\newcommand{\IR}{{\mathbb{R}}}
\newcommand{\IN}{{\mathbb{N}}}

\newcommand{\IS}{{\mathbb{S}}}

\newcommand{\cB}{\mathcal{B}}

\newcommand{\cT}{\mathcal{T}}

\usepackage{xcolor}

\usepackage{slashed} 
\newcommand{\D}{\slashed{D}}

\DeclareMathOperator{\spec}{spec}

\DeclareMathOperator{\dom}{dom}

\DeclareMathOperator{\id}{id}

\DeclareMathOperator{\vol}{vol}

\usepackage[backend=bibtex, doi=false, isbn=false, url=false]{biblatex}
\usepackage{xcolor}

\title[On an iterative scheme for the spinorial Yamabe Equation]{On an Iterative Scheme for the Spinorial Yamabe Equation
	on Manifolds with Boundary}
\author{Eric Trébuchon}
\date{}
\thanks{Mathematisches Institut, Universität Freiburg, \texttt{eric.trebuchon@math.uni-freiburg.de}}

\begin{document}
	
\begin{abstract}
	We develop a general fixed-point argument for establishing the existence of solutions to non-linear boundary value problems involving first-order elliptic operators on vector bundles, with a primary application to the spinorial Yamabe and Brezis--Nirenberg equations on compact manifolds with boundary. Motivated by the iterative approach of Rosenberg and Xu \cite{Rosenberg21}, we extend and improve their techniques, applied to the spinorial setting. We prove the existence of an $H^1$-solution under an explicit smallness condition involving the boundary datum and the inverse operator, and the uniqueness of the solution in a ball in a suitable Lebesgue space.
	The iteration converges geometrically in $L^q$ and $H^1$, and the solution depends continuously on the harmonic representative of the boundary datum. The smallness condition is homogeneous under amplitude scaling and simultaneous scaling of the operator and the non-linear coupling.
	Furthermore, we provide $C^{1,\alpha}$-regularity in the interior and smoothness away from the zero set by bootstrapping, with corresponding conclusions up to the boundary under appropriate elliptic regularity assumptions and sufficient regularity of the boundary data.
	Our approach is not limited to the Dirac operator and provides a versatile tool for a wide class of non-linear problems.
	Finally, we explain its application to families of small domains under uniform inverse estimates after rescaling to a fixed reference domain.
\end{abstract}

\maketitle
\tableofcontents

\section{Introduction}

The study of spin conformal invariants has played a central role in differential geometry, with foundational contributions by Hijazi, Lott, and others \cite{Hijazi86,Lott86,Baer92,Ammann03hab}, who introduced such an invariant via the smallest positive eigenvalue of the Dirac operator on closed manifolds. This concept was later extended to manifolds with boundary by Raulot \cite{Raulot08}, who introduced a conformal invariant associated to a fixed spin structure $\sigma$ on a compact Riemannian manifold $(\Omega, g)$ with boundary, under the chiral boundary condition $B^\pm\psi=0$, defined by
\begin{equation}
	\lambda_{\text{min}} (\Omega, \partial \Omega) \coloneq \inf_{g' \in [g]} |\lambda_1 (\D^{g'}_{B^\pm})| \, \vol(\Omega, g')^{\frac{1}{n}},
	\label{def:lambda_min}
\end{equation}
where $[g]$ represents the conformal class of the Riemannian metric $g$ and $\lambda_1(\D^g_{B^\pm})$ denotes the eigenvalue of smallest modulus of the Dirac operator $\D_g$ under the chiral boundary condition $B^\pm$.
Here we work in the invertible setting. Statements involving chirality are understood on a bundle carrying a chirality operator; in odd dimension this may require passing to a doubled spinor bundle.
As $\D^g_{B^+}$ and $\D^g_{B^-}$ have the same spectrum up to sign, the choice of sign is irrelevant for $\lambda_{\text{min}} (\Omega, \partial \Omega)$.

This article mainly addresses the spinorial Yamabe equation on Riemannian spin manifolds with boundary. This equation is associated with the variational characterization of the above-mentioned conformal invariant.
Specifically, given a smooth compact Riemannian manifold $(\Omega, g)$ with smooth boundary $\partial \Omega$, dimension $n\geq2$, a fixed spin structure $\sigma$ and its corresponding spinor bundle $\Sigma \Omega$, as well as a chiral boundary condition $B^\pm$, the Euler--Lagrange equation turns out to be equivalent with the transformations and normalization explained in Appendix~\ref{appendix} to
\begin{align}
	\begin{split}
		\D u &= \lambda |u|^{p-2} u \quad \text{in } \Omega, \label{def:pde_homogenous}\\
		B^\pm u &= 0\quad \quad \,\,\,\,\,\text{on } \partial\Omega.
	\end{split}
\end{align}
Here $p = \frac{2n}{n-1}$ is the critical Sobolev exponent for the embedding $H^{1/2}\hookrightarrow L^p$ and $\lambda\in \IR$.  It is also the right exponent for the problem to be conformally invariant.

In our article, however, we will analyse this problem for a wider range of boundary operators $P$ and with non-homogeneous boundary condition, meaning we will be interested in showing the existence of solutions to 
\begin{align}
	\begin{split}
		\D u &= \lambda |u|^{p-2} u \quad \text{in } \Omega, \label{def:pde_inhomogenous}\\
		Pu &= P\psi_0\quad \quad \,\,\,\,\,\text{on } \partial\Omega,
	\end{split}
\end{align}
where $\psi_0\in H^1(\Sigma\Omega)$ is a prescribed section of the spinor bundle and $P$ is a bounded linear boundary operator such that
	\[
	\D_P\coloneq\D\big|_{\dom(\D_P)},
	\qquad
	\dom(\D_P)=\{v\in H^1(\Sigma\Omega):Pv=0\},
	\]
	is self-adjoint, regular and invertible.
	
Our methods also work for a wider range of values for $p$, see Remark~\ref{remark: generality of method}. The argument extends beyond the critical exponent and
the Dirac operator. Its essential ingredients are a
bounded inverse from $L^2$ to $H^1$, a suitable Sobolev
embedding, and growth and Lipschitz estimates for the
nonlinearity. We discuss other
exponents, variable coefficients, and more general
nonlinearities in
Remark~\ref{remark: generality of method}.

The classical analogue, the Yamabe problem on manifolds with boundary, has been studied extensively using a variety of analytical methods \cite{Brendle14,Escobar92_1,Escobar92_2,Fernando05}.
\newline
In the following, we first want to state some results for the conformal spin invariant on closed manifolds, which is defined as
\[
\lambda_{\text{min}}^+(\Omega,[g])
=
\inf_{g'\in[g]}
\la_1^+(\D^{g'})
\t{Vol}(\Omega,g')^{\f1n},
\]
where $\la_1^+$ denotes the first positive eigenvalue.
In \cite{Hijazi91}, Hijazi derives a conformal lower bound for the eigenvalues $\lambda$ of the Dirac operator involving the Yamabe invariant $\mu(\Omega,[g])$, namely if $n\geq3$, then
\[
\lambda^2\t{Vol}(\Omega,g)^{\f2n}
\geq \f{n}{4(n-1)}\mu(\Omega,[g]).
\]
The Yamabe invariant $\mu(\Omega,[g])$ has been introduced in \cite{Yamabe60} in order to solve the Yamabe problem of finding a conformal metric with constant scalar curvature, although this paper contained a flaw pointed out by Trudinger, which was then corrected by Trudinger, Aubin and later in full generality by Schoen, resulting in a complete affirmative resolution of the Yamabe problem in the year 1984.

Lott showed in \cite{Lott86} that the conformal invariant is positive if there is some element $g'\in[g]$ such that $\D^{g'}$ is invertible. This result was improved in \cite{Ammann03}, where Ammann showed that the assumption can be dropped.
In the same paper, Ammann showed that
\[
\lambda_{\text{min}}^+(\Omega,[g])
\leq\lambda_{\text{min}}^+(\IS^n)
=\f n2\omega_n^{\f1n}
\]
with the standard metric and spin structure on $\IS^n$, where $\omega_n$ is the $n$-measure of $\IS^n$. In \cite{Ammann03hab}, he showed that the strict inequality
\[
\lambda_{\text{min}}^+(\Omega,[g])
<
\lambda_{\text{min}}^+(\IS^n)
\]
implies the existence of a solution to the spinorial Yamabe equation. This is very reminiscent of the result by Aubin \cite{Aubin76} in resolving the classical Yamabe problem. Several other spinorial Yamabe-type equations can be found, e.g., in \cite{Isobe11_1,Isobe11_2,Raulot09,Borrelli19}.

Solutions to the spinorial Yamabe equation are connected to several geometric problems. For example, in dimension $n=2$, the spinorial representation relates solutions to conformal immersions into $\IR^3$ with constant mean curvature \cite{Ammann03}.
The spinorial Yamabe equation is also related to non-linear Dirac equations occurring in the theory of Dirac-harmonic maps \cite{Wang06,Wang08}. Dirac operators with boundary conditions also occur in
non-linear effective models of graphene \cite{AkhmerovBeenakker08}.
These provide further motivation for studying Dirac
boundary value problems.
\newline
Finally arriving back at our setting with both the boundary and the spinorial version, the spin conformal invariant~\eqref{def:lambda_min} introduced in \cite{Raulot08} satisfies several analogous results:
\begin{enumerate}
	\item The Hijazi inequality proved in \cite{Raulot06} yields a comparison between the chiral invariant $\la_\t{min}(\Omega,\partial\Omega)$ and the Yamabe invariant $\mu(\Omega,\partial\Omega)$ (see \cite{Escobar92_2}). Explicitly, if $n\geq3$, we have
	\[
	\la_\t{min}(\Omega,\partial\Omega)^2
	\geq\f{n}{4(n-1)}\mu(\Omega,\partial\Omega).
	\]
	\item In \cite{Raulot08}, a generalization of Ammann's result, which holds for closed spin manifolds, is shown:
	\[
	\la_\t{min}(\Omega,\partial\Omega)
	\leq
	\lambda_{\text{min}}(\IS^n_+,\partial\IS^n_+)
	=
	\f n2\left(\f{\om_n}{2}\right)^{\f1n},
	\]
	where $\om_n$ again denotes the $n$-dimensional volume of the sphere $\IS^n$.
	\item For certain classes of manifolds, the last inequality is shown to be strict \cite{Raulot11}. There are similar results in the context of closed manifolds.
	In the classical Yamabe problem with boundary, strict comparison with the model invariant plays a corresponding role in existence results; see \cite{Escobar92_2}.
	\item For spinorial Yamabe-type equations on manifolds with boundary, there are some known results, e.g. \cite{Ding18,Wen24}. The latter employs variational methods to find critical points of the associated functional.
	In contrast, the present work focuses on a local small-data boundary value problem. By employing an iterative scheme, we provide existence, uniqueness in a specified ball, continuous dependence, and a convergent sequence of approximations. These conclusions do not assert global uniqueness and do not show existence for homogenous boundary datum.
\end{enumerate}

We now return to our problem of solving~\eqref{def:pde_inhomogenous}. In \cite{Rosenberg21}, the authors gave a new approach for solving the classical Yamabe equation with Dirichlet boundary conditions under some smallness assumptions on the underlying Riemannian domain using an iterative scheme. Iterative schemes have been used for non-linear PDEs in the past \cite{Moser61,Moser66}.

The main objective of this paper is to develop a general fixed-point/iterative framework for solving non-linear boundary value problems involving first-order elliptic operators on vector bundles. While we apply this framework to the spinorial Yamabe equation with boundary, the resulting methods turn out to be a versatile toolbox for a wide class of non-linear partial differential equations. We discuss in Remark~\ref{remark: generality of method} to what degree we can change not only the differential operator and the boundary condition, but also the critical exponent $p=\f{2n}{n-1}$, including the range
	\[
	2<p\leq\frac{2n-2}{n-2}\quad(n>2)
	\]
	and every finite $p>2$ when $n=2$. The case $p=2$ is linear and is treated separately.

A further application is the spinorial
Brezis--Nirenberg-type boundary value problem
\[
\D u=a u+\lambda|u|^{p-2}u,
\qquad
Pu=P\psi_0,
\qquad
p=\frac{2n}{n-1}.
\]
For $a\in\mathbb R\setminus\operatorname{spec}(\D_P)$,
the linear term can be absorbed into the operator:
one replaces $\D_P$ by $\D_P-a$ and uses the
$(\D-a)$-harmonic representative of the boundary datum.
The same fixed-point argument then yields existence,
uniqueness in a prescribed ball, continuous dependence
on the boundary data, and geometric convergence under
the corresponding smallness condition.
The regularity conclusions also remain valid; see
Corollary~\ref{cor:linear-perturbation}.
This extension concerns small inhomogeneous boundary
data and does not assert the existence of nonzero
solutions with homogeneous boundary data.

While the spinorial Yamabe equation is naturally associated with a variational problem, our purpose is not to construct a global variational minimizer. Instead, the iterative scheme presented in this paper relies on smallness assumptions involving the prescribed boundary datum and the inverse of the operator.

In the Yamabe case, working with the conformal Laplace operator, many estimates are carried out on the corresponding quadratic forms using $H^1$-regularity. Although $H^{1/2}$ is the natural variational space,
we perform the contraction directly in $L^q$.
This avoids the use of fractional Sobolev spaces
in the fixed-point argument, while $H^1$ regularity
follows from the mapping properties of the inverse.
For the iteration, we use an exponent adapted to the
	nonlinearity. Set
	\[
	F(v)=|v|^{p-2}v,
	\qquad
	q=2(p-1)=\frac{2(n+1)}{n-1}.
	\]
	Then
	\[
	F:L^q(\Sigma\Omega)\longrightarrow L^2(\Sigma\Omega),
	\qquad
	\|F(v)\|_{L^2}=\|v\|_{L^q}^{p-1}.
	\]
	The inverse $\D_P^{-1}$ maps $L^2$ boundedly into $H^1$,
	and the Sobolev embedding $H^1\hookrightarrow L^q$
	therefore makes the fixed-point map well defined on $L^q$.
	
	To estimate the inverse norm quantitatively, we use
	the Gagliardo--Nirenberg inequality
	\[
	\|v\|_{L^q}
	\leq C_{\mathrm{GN}}
	\|v\|_{L^2}^{1-s}\|v\|_{H^1}^{s},
	\qquad
	s=\frac n{n+1},
	\]
	where $1/q=1/2-s/n$.

To show our upcoming results, we will need that $\D_P$ is a \textit{regular} operator, by which we mean that there exists a $c_1>0$ such that for all $v\in\dom(\D_P)$

	\begin{equation*}
		\|v\|_{H^1}^2
		\leq c_1\left(\|v\|_{L^2}^2+\|\D_Pv\|_{L^2}^2\right).
	\end{equation*}

	Since the boundary value problem only depends on $P\psi_0$, we first replace $\psi_0$ by its harmonic representative
	\[
	\psi=\psi_0-\D_P^{-1}(\D\psi_0).
	\]
	It satisfies $\D\psi=0$, preserves $P\psi=P\psi_0$, and is unique with this property by invertibility of $\D_P$.
	
	The constant governing the fixed-point argument is $
	A_q\coloneq\|\D_P^{-1}\|_{\mathcal L(L^2,L^q)}.$ It can be estimated by geometric and spectral datum.

With all this, our theorem can be stated in its simplest form as:	

\begin{theorem}[Existence, uniqueness, and stability]
	\label{thm:main_stability}
	Let $\D_P$ be a self-adjoint, invertible, and regular
	realization of $\D$ on a smooth compact Riemannian
	spin manifold of dimension $n\geq2$ with boundary.
	Set
	\[
	p=\frac{2n}{n-1},
	\qquad q=2(p-1),
	\qquad
	A_q=\|\D_P^{-1}\|_{\mathcal L(L^2,L^q)}.
	\]
	Let $\psi\in H^1(\Sigma\Omega)$ satisfy $\D\psi=0$.
	If
	\[
	|\lambda|A_q\|\psi\|_{L^q}^{p-2}
	<
	\frac{(p-2)^{p-2}}{(p-1)^{p-1}},
	\]
	then the boundary value problem
	\[
	\D u=\lambda|u|^{p-2}u,
	\qquad Pu=P\psi
	\]
	has an $H^1$-solution, unique among those satisfying
	\[
	\|u-\psi\|_{L^q}
	\leq\frac{\|\psi\|_{L^q}}{p-2}.
	\]
	
	The solution depends continuously on the harmonic
	representative of the boundary datum.
	More precisely, for every $M>0$ satisfying
	\[
	|\lambda|A_qM^{p-2}
	<
	\frac{(p-2)^{p-2}}{(p-1)^{p-1}},
	\]
	there is a constant $C_M>0$ such that the solutions
	corresponding to any two harmonic representatives
	$\psi_1,\psi_2$ with $\|\psi_j\|_{L^q}\leq M$ satisfy
	\[
	\|u_1-u_2\|_{L^q}
	\leq C_M\|\psi_1-\psi_2\|_{L^q}.
	\]
	Convergence of the representatives in $H^1$ also
	implies convergence of the corresponding solutions
	in $H^1$.
	
	For every $u_0\in L^q$ in the above $L^q$-ball around $\psi$, the iteration
	\[
	u_{k+1}
	=
	\psi+\lambda\D_P^{-1}
	\bigl(|u_k|^{p-2}u_k\bigr)
	\]
	converges geometrically to $u$ in $L^q$ and $H^1$.
\end{theorem}

Quantitative convergence and stability
are given in Theorems~\ref{prop:quantitative}
and~\ref{prop:stability}.

Inspired by the methods in \cite{Rosenberg21} for the classical Yamabe equation, we use similar ideas, though we present the proof as a fixed-point argument instead, which gives uniqueness in the contraction ball, continuous dependence on the boundary datum, and quantitative convergence. Whenever the corresponding contraction estimates are available in the classical setting, the same fixed-point principle gives analogous local uniqueness and stability statements there. Another major difference is our choice of not working with the associated energy-space $H^{\f12}$, and circumvent the associated prerequisites and technicalities.

Let $\mu=\operatorname{dist}(0,\operatorname{spec}(\D_P))>0$.
	The spectral gap enters through an estimate for the inverse. In particular, the Gagliardo--Nirenberg inequality and regularity give
	\[
	A_q
	\leq C_{\mathrm{GN}}c_1^{s/2}
	\frac{(1+\mu^2)^{s/2}}{\mu}.
	\]
	If the remaining constants are controlled, the reciprocal of this bound behaves like $\mu$ for small $\mu$ and like $\mu^{1/(n+1)}$ for large $\mu$.
	At the same time, inserting a normalized eigenspinor whose eigenvalue has modulus $\mu$, and using Hölder gives
	\[
		A_q\geq \frac{\vol_g(\Omega)^{1/q-1/2}}{\mu}
	=\frac{\vol_g(\Omega)^{-1/(n+1)}}{\mu},
	\]
	meaning the lower bound shows the limitation of this constant; the reciprocal behaves indeed like $\mu$ for small $\mu$.

We adapt these techniques to the spinorial setting, but the Dirac operator does not have an analogous domain-monotonicity principle as the Laplacian: the Dirichlet Rayleigh quotient and continuation by zero give domain-monotonicity for the Dirichlet Laplacian. There is no general corresponding principle for the Dirac boundary realizations that we can use here. This is part of why we formulate the conditions in terms of inverse and regularity estimates.

The argument also applies to suitable nonsmooth domains,
see Remark~\ref{rem:not smooth}. For our iterative scheme,
the relevant Sobolev embedding and a bounded inverse
$L^2\to H^1$ suffice. For suitable planar corner domains,
self-adjointness and regularity of Dirac operators with local boundary conditions
are studied in \cite{VanDenBosch21}.
Invertibility and the required embedding estimates
must also be verified.

It is important to distinguish our result from
existence of nonzero solutions with homogeneous
boundary data.
If $P\psi_0=0$, its harmonic representative is zero,
and the small-data contraction constructs only the
zero solution.
Nevertheless, invertibility gives a quantitative
gap between zero and every nontrivial
$H^1$-solution of the homogeneous problem:
\[
\|u\|_{L^q}
\geq
\bigl(|\lambda|A_q\bigr)^{-1/(p-2)},
\qquad q=2(p-1).
\]
Indeed, this follows directly from the fixed-point
identity and
$\|u\|_{L^q}\leq|\lambda|A_q\|u\|_{L^q}^{p-1}$.
Thus the method also yields an explicit $L^q$-norm
gap for nontrivial solutions; see
Corollary~\ref{cor:Lq-gap} and Remark~\ref{rem:homogeneous-gap-Brezis}.
For nonzero prescribed boundary data, the solution
constructed by our existence theorem is necessarily
nonzero. These statements hold for both the spinorial Yamabe problem, as well as the spinorial Brezis--Nirenberg-type problem.

In Section~\ref{Regularity via Bootstrapping}, we enhance
the regularity by iterating Sobolev embeddings and
applying Schauder estimates. Interior regularity does
not require additional assumptions on the boundary
operator. For regularity up to the boundary, however,
we require stronger elliptic estimates.

More precisely, for global Sobolev regularity we assume
the following implication for every $1<a<\infty$:
if $v$ satisfies the homogeneous boundary condition
and $v,\D v\in L^a$ in the weak sense, then
$v\in W^{1,a}$ and
\begin{equation*}
	\|v\|_{W^{1,a}}
	\leq C_a\left(
	\|v\|_{L^a}+\|\D v\|_{L^a}
	\right).
\end{equation*}
For smooth local elliptic boundary problems satisfying
the Shapiro--Lopatinski condition, such regularity
statements are available from elliptic boundary theory.
Boundary Schauder estimates are an additional requirement
for the corresponding conclusions in H\"older spaces.

\begin{theorem}[Regularity]
	\label{second main theorem}
	Let $(\Omega,g)$ be a smooth compact Riemannian spin manifold of dimension $n\geq2$ with smooth boundary.
	Suppose that $u\in H^1(\Sigma\Omega)$ satisfies
	\[
	\D u=\lambda|u|^{p-2}u,
	\qquad
	p=\frac{2n}{n-1}.
	\]
	Then $u\in C_{\mathrm{loc}}^{1,\alpha}\bigl(\Sigma(\operatorname{int}\Omega)\bigr)$ for every $0<\alpha<1$, and $u$ is smooth in the interior away from its zero set.
	
	Furthermore, assume that $Pu=P\psi$ and that the boundary operator $P$ satisfies the Shapiro--Lopatinski condition relative to $\D$. Then:
	\begin{enumerate}
		\item If the prescribed boundary datum $\psi$ admits an extension in $\bigcap_{1<a<\infty} W^{1,a}(\Sigma\Omega)$, then $u\in \bigcap_{1<a<\infty} W^{1,a}(\Sigma\Omega)$.
		\item If $P$ additionally admits boundary Schauder estimates and $\psi$ admits a $C^{1,\alpha}$ extension, then $u\in C^{1,\alpha}(\Sigma\Omega)$.
		\item If the datum $\psi$ is smooth and $P$ admits boundary Schauder regularity of all orders, then $u$ is smooth away from its zero set up to the boundary.
	\end{enumerate}
\end{theorem}

We note that besides the chiral boundary condition,
other local boundary conditions giving rise to
self-adjoint and regular Dirac operators exist;
these are studied in
\cite{Grosse-Uribe-vanDenBosch}.

\paragraph{\textbf{AI Disclosure.}}
During the preparation of this manuscript, the author used large language models (Qwen 3.5 and OpenAI GPT-6 Astra) to check author-written proofs for typographical errors, notational inconsistencies, and mathematical mistakes, and to improve the write-up of the main result. These tools also helped identify consequences of the fixed-point argument, refine the condition: contractions using a fractional estimate, and introduce the constant \(A_q\) to express the conditions homogeneously under scaling of \(\D_P\).  All other mathematical content, results, and arguments are the author's own, and the author carefully checked all revisions and takes full responsibility for the mathematical content.
\newpage

	\section{Prerequisites}
	
	\subsection{Sobolev Spaces on Vector Bundles}
	
	Let $E$ be a smooth Hermitian vector bundle of rank $r$
	over an $n$-dimensional compact Riemannian manifold
	$(\Omega,g)$ with smooth boundary, equipped with a smooth
	Hermitian connection $\nabla^E$.
	
	For $1\leq a<\infty$, the intrinsic Lebesgue norm is
	\[
	\|u\|_{L^a(E)}
	=
	\left(\int_\Omega |u|^a\,d\vol_g\right)^{1/a}.
	\]
	For $a=\infty$, we use the essential supremum norm.
	In particular, $L^2(E)$ carries the Hermitian product
	\[
	\langle u,v\rangle_{L^2}
	=
	\int_\Omega\langle u,v\rangle\,d\vol_g.
	\]
	
	For an integer $m\geq1$ and $1<a<\infty$, the Sobolev
	space $W^{m,a}(E)$ consists of sections whose weak
	covariant derivatives up to order $m$ belong to $L^a$.
	We use the norm
	\[
	\|u\|_{W^{m,a}(E)}^a
	=
	\sum_{j=0}^m
	\int_\Omega |(\nabla^E)^j u|^a\,d\vol_g.
	\]
	The space $W^{m,\infty}(E)$ is defined analogously using
	essential supremum norms. We write
	\[
	H^1(E)=W^{1,2}(E),
	\qquad
	\|u\|_{H^1(E)}^2
	=
	\|u\|_{L^2(E)}^2+\|\nabla^E u\|_{L^2(E)}^2.
	\]
	
	These definitions agree with the usual local definitions
	of integer-order Sobolev spaces, with equivalent norms;
	see \cite{Große13}.
	H\"older spaces $C^{m,\alpha}(E)$ are defined using local
	representations of sections, but will only appear in Theorem~\ref{second main theorem}.
	
	Since $\Omega$ has finite volume, H\"older's inequality gives
	\[
	\|u\|_{L^a(E)}
	\leq
	\vol_g(\Omega)^{1/a-1/b}\|u\|_{L^b(E)},
	\qquad 1\leq a<b\leq\infty.
	\]
	In particular, $L^b(E)\hookrightarrow L^a(E)$ and
	$W^{1,b}(E)\hookrightarrow W^{1,a}(E)$ continuously
	whenever $1<a<b\leq\infty$.
	
	\subsection{Sobolev Embeddings and the
		Gagliardo--Nirenberg Inequality}
	
	We use the following standard embeddings, which extend
	to vector bundles by localization and trivialization of the vector bundles
	\cite{Aubin82,Große13}.
	
	\begin{theorem}[Sobolev embeddings]
		\label{sobolev embeddings}
		The following embeddings are continuous:
		\begin{enumerate}
			\item[(i)]
			$W^{1,a}(E)\hookrightarrow L^{na/(n-a)}(E)$
			for $1<a<n$, and
			$W^{1,n}(E)\hookrightarrow L^b(E)$
			for every finite $b\geq1$.
			\item[(ii)]
			$W^{1,a}(E)\hookrightarrow C^{0,\alpha}(E)$
			for $a>n$ and $0<\alpha\leq1-n/a$.
		\end{enumerate}
		Moreover, $H^1(E)\hookrightarrow L^2(E)$ is compact.
	\end{theorem}
	
	In particular,
	\[
	H^1(E)\hookrightarrow L^q(E)
	\]
	for $2\leq q\leq 2n/(n-2)$ when $n>2$, and for every
	finite $q\geq2$ when $n=2$.
	
	The following will be necessary to prove the Gagliardo--Nirenberg inequality:
	\begin{lemma}[First-order Kato inequality]
		\label{lem:kato}
		Let $E$ be a smooth Hermitian vector bundle over
		$(\Omega,g)$, equipped with a metric connection
		$\nabla^E$. If $u\in H^1(E)$, then
		$|u|\in H^1(\Omega)$ and
		\[
		|d|u||\leq|\nabla^E u|
		\qquad\text{almost everywhere}.
		\]
		In particular,
		\[
		\||u|\|_{L^2(\Omega)}=\|u\|_{L^2(E)},
		\qquad
		\||u|\|_{H^1(\Omega)}\leq\|u\|_{H^1(E)}.
		\]
	\end{lemma}
	
	\begin{proof}
		For $\varepsilon>0$, define
		\[
		f_\varepsilon=(|u|^2+\varepsilon^2)^{1/2}.
		\]
		The Sobolev chain rule, applied in smooth local
		orthonormal frames, gives $f_\varepsilon\in H^1(\Omega)$.
		Metric compatibility implies
		\[
		df_\varepsilon(X)
		=
		\frac{\operatorname{Re}
			\langle\nabla_X^E u,u\rangle}
		{(|u|^2+\varepsilon^2)^{1/2}}
		\]
		almost everywhere.
		Consequently, by Cauchy--Schwarz,
		\[
		|df_\varepsilon|^2
		\leq
		\frac{|u|^2}{|u|^2+\varepsilon^2}
		|\nabla^E u|^2
		\leq|\nabla^E u|^2.
		\]
		
		Since
		\[
		0\leq f_\varepsilon-|u|\leq\varepsilon
		\]
		and $\Omega$ has finite volume, we have
		$f_\varepsilon\to|u|$ in $L^2(\Omega)$.
		The preceding derivative estimate shows that
		$f_\varepsilon$ is bounded in $H^1(\Omega)$ as
		$\varepsilon\downarrow0$.
		Thus, along a subsequence, it converges weakly
		in $H^1$ to $|u|$, proving that $|u|\in H^1(\Omega)$.
		
		For every measurable subset $A\subset\Omega$,
		weak lower semicontinuity gives
		\[
		\int_A |d|u||^2\,d\vol_g
		\leq
		\liminf_{\varepsilon\downarrow0}
		\int_A |df_\varepsilon|^2\,d\vol_g
		\leq
		\int_A|\nabla^E u|^2\,d\vol_g,
		\]
		where the limit inferior is taken along the chosen
		subsequence.
		Hence $|d|u||\leq|\nabla^E u|$ almost everywhere.
		The norm statements follow immediately.
	\end{proof}
	
	The next estimate will be used to bound
	$A_q \coloneq \|\D_P^{-1}\|_{\mathcal L(L^2,L^q)}$.
	
	\begin{lemma}[Gagliardo--Nirenberg inequality]
		\label{Fractional Gagliardo--Nirenberg}
		Let $0<s<1$ and
		\[
		\frac1q=\frac12-\frac{s}{n}>0.
		\]
		There exists $C_{\mathrm{GN}}>0$ such that for all $u\in H^1(E)$,
		\begin{equation}
			\|u\|_{L^q(E)}
			\leq C_{\mathrm{GN}}
			\|u\|_{L^2(E)}^{1-s}\|u\|_{H^1(E)}^s.
			\label{GN inequality}
		\end{equation}
	\end{lemma}
	
	\begin{proof}
		The classical scalar Gagliardo--Nirenberg inequality on
		a compact manifold with smooth boundary gives
		\[
		\|f\|_{L^q}
		\leq C
		\|f\|_{L^2}^{1-s}\|f\|_{H^1}^s.
		\]
		It follows from the Euclidean inequality by localization
		and extension across the boundary; see
		\cite{Aubin82,Hajaiej11}.
		
		For $u\in H^1(E)$, using Lemma~\ref{lem:kato}, the scalar function $|u|$ belongs
		to $H^1(\Omega)$ and satisfies the inequality $|d|u||\leq|\nabla^E u|$ a.e. on $\Omega$ and hence applying the Gagliardo--Nirenberg inequality to $|u|$ and using
		\[
		\||u|\|_{L^2}=\|u\|_{L^2(E)},\qquad \||u|\|_{L^q}=\|u\|_{L^q(E)},
		\qquad
		\||u|\|_{H^1}\leq\|u\|_{H^1(E)}
		\]
		proves
		\eqref{GN inequality}.
	\end{proof}
	
	The constant $C_{\mathrm{GN}}$ depends on the geometry,
	the exponent and the chosen norms. Uniform estimates
	for families of domains require separate control of
	these constants.
	
	\subsection{Trace Operator}
	
	For $1<a<\infty$, the trace extends to a bounded
	surjective operator
	\[
	\operatorname{Tr}_{\partial\Omega}:
	W^{1,a}(E)\longrightarrow
	W^{1-1/a,a}(E|_{\partial\Omega}),
	\]
	which admits a bounded linear right inverse
	\cite{Schwarz95}. For continuous sections in
	$W^{1,a}(E)$, it agrees with restriction to the boundary.
	
	In particular, $H^1(E)$ has a bounded trace in
	$H^{1/2}(E|_{\partial\Omega})$.
	
	Boundary conditions are imposed on the trace of a
	Sobolev section. In particular, a smooth field of
	projections $	\pi:E|_{\partial\Omega}\longrightarrow E|_{\partial\Omega}$ 
	defines a local boundary operator by $Pu=\pi\operatorname{Tr}_{\partial\Omega}u. $
	The chiral boundary operators introduced below are
	of this form.
	
	\subsection{The Spinor Bundle}
	
	Let $(\Omega,g,\sigma)$ be an $n$-dimensional compact
	Riemannian spin manifold.
	The spinor bundle is denoted by
	\[
	\Sigma\Omega=\Sigma(\Omega,g,\sigma)
	\]
	and carries a natural Clifford multiplication $\cdot$,
	a natural Hermitian metric $\langle\cdot,\cdot\rangle$
	and a metric connection $\nabla^S$. These fulfil the
	properties of a Dirac bundle as introduced in
	\cite{Lawson89}.
	We choose our Clifford relation to be with the
	$+$-convention:
	\[
	v\cdot w\cdot\psi+w\cdot v\cdot\psi
	=
	2g(v,w)\psi.
	\]
	Clifford multiplication by a tangent vector is
	then Hermitian.
	With these notions, we may define the Dirac operator
	by acting as
	\[
	\D=-i\sum_{k=1}^n e_k\cdot\nabla^S_{e_k}
	\]
	for local orthonormal frames $(e_k)_{k=1,\dots,n}$.
	This definition is independent of the local
	orthonormal frame chosen.
	
	We consider boundary realizations of $\D$ on $L^2$
	with domains of the form
	\[
	\dom(\D_P)
	=
	\{\varphi\in H^1(\Sigma\Omega):P\varphi=0\}.
	\]
	The assumptions on $P$ and the corresponding
	regularity estimate are stated below.
	
	Denote the rank of the spinor bundle by
	$r\coloneq\operatorname{rank}(\Sigma\Omega)$, which
	always equals $2^{\lfloor n/2\rfloor}$ for the classical
	spinor bundle of a spin manifold.
	
	The pointwise estimate
	\[
	|\D\varphi|
	\leq\sqrt n\,|\nabla^S\varphi|
	\]
	shows that
	\[
	\D:W^{1,a}(\Sigma\Omega)\longrightarrow L^a(\Sigma\Omega)
	\]
	is bounded for $1<a\leq\infty$.
	
	\subsection{Chiral Boundary Condition}
	We will define the chiral boundary condition following
	\cite{Raulot08,Baer12}, with the Clifford convention
	specified above.
	
	We retain the conventions of the preceding subsection
	and write
	\[
	\gamma(X)\varphi \coloneq X\cdot\varphi
	\]
	for Clifford multiplication.
	Let $\nu$ denote the outward unit normal along
	$\partial\Omega$.
	
	A chirality operator is a parallel unitary bundle
	endomorphism $\beta$ satisfying
	\[
	\beta^2=\id,
	\qquad
	\beta\gamma(X)=-\gamma(X)\beta
	\]
	for every tangent vector $X$.
	Parallelness means that
	\[
	\nabla_X^S(\beta\varphi)
	=
	\beta\nabla_X^S\varphi.
	\]

	Such an operator always exists if $\Omega$ is an
	even-dimensional spin manifold and $\Sigma$ the
	standard spinor bundle. If $\Omega$ were
	odd-dimensional, one can pass to the double of the
	spinor bundle
	$\widetilde\Sigma\coloneq\Sigma\oplus\Sigma$
	to ensure the existence of such a $\beta$:
	\begin{enumerate}
		\item
		For even $n$, given a spin manifold $\Omega$ with
		spinor bundle $\Sigma$, we can choose
		\[
		\beta=i^{n/2}\gamma(e_1)\cdots\gamma(e_n)
		\]
		as a chirality operator for a positively oriented
		orthonormal frame.
		
		\item
		For odd $n$, we generally do not have a chirality
		operator on the standard spinor bundle $\Sigma$
		over $\Omega$, but the following work-around suffices.
		Consider the spinor bundle $\Sigma$ of complex rank
		$r=2^{\lfloor n/2\rfloor}$. To ensure the existence of
		a chirality operator $\beta$, we define the augmented
		spinor bundle
		$\widetilde\Sigma\coloneq\Sigma\oplus\Sigma$
		with the product connection and choose an extended
		Clifford multiplication by
		\[
		\widetilde\gamma(X)
		=
		\begin{pmatrix}
			0&\gamma(X)\\
			\gamma(X)&0
		\end{pmatrix}.
		\]
		This construction allows us to define a global
		chirality operator
		\[
		\beta=
		\begin{pmatrix}
			\id&0\\
			0&-\id
		\end{pmatrix},
		\]
		and the $\pm1$ eigenspaces split the bundle into the
		same splitting
		$\widetilde\Sigma=\Sigma\oplus\Sigma$.
	\end{enumerate}
	
	The chiral bag boundary condition is defined by
	projecting spinors onto the eigensubbundles $V_\pm$
	of the boundary spinor bundle, corresponding to the
	eigenvalues $\pm i$ of the unitary map
	$\gamma(\nu)\beta$, which squares by our convention
	to $-\id$.
	The projection operators are
	\[
	P_{\t{CHI}}^\pm:
	H^1(\Sigma\Omega)\to H^{1/2}(V_\pm),
	\qquad
	\varphi\mapsto
	\frac12(\id\mp i\gamma(\nu)\beta)
	\varphi|_{\partial\Omega}.
	\]
	These operators define two boundary conditions
	$P_{\t{CHI}}^\pm(\varphi)=0$. These are the boundary operators denoted by
	$B^\pm$ in the introduction: $
	B^\pm=P_{\mathrm{CHI}}^\pm. $
	
	\begin{remark}[Terminology]
		Terminology varies between the geometric and physical
		literature. Throughout this paper, the chiral boundary
		condition means $P_{\mathrm{CHI}}^\pm$ as defined above.
		Related (often the same) boundary conditions occur under the names
		``MIT bag'' and ``infinite mass'', but their formulas
		depend on the Clifford representation and conventions:			
		In our convention for example, the Riemannian MIT boundary
		projectors are $
		Q^\pm=\frac12\bigl(\id\pm\gamma(\nu)\bigr).$
		They differ from $P_{\mathrm{CHI}}^\pm$ and define
		non-self-adjoint realizations of the Dirac operator. 
		
		For physical background and infinite-mass limits, see
		\cite{Bogolioubov68,Chodos74,Johnson75,Arrizabalaga19,
			BerryMondragon87,AkhmerovBeenakker08}.
		For local self-adjoint boundary conditions and their
		regularity, see
		\cite{Grosse-Uribe-vanDenBosch}.
		Non-local boundary conditions are also treated in
		\cite{Baer12}.
	\end{remark}

	\subsection{Regularity and Spectral Properties}
	
	The following arguments apply to any self-adjoint
	realization $\D_P$ with domain
	\[
	\dom(\D_P)
	=
	\{u\in H^1(\Sigma\Omega):Pu=0\}
	\]
	satisfying the regularity estimate
	\begin{equation}
		\|u\|_{H^1}^2
		\leq
		c_1\bigl(\|u\|_{L^2}^2+\|\D_Pu\|_{L^2}^2\bigr),
		\qquad u\in\dom(\D_P).
		\label{regularity spectral}
	\end{equation}
	Here $P$ is a bounded linear boundary operator on
	$H^1(\Sigma\Omega)$.
	
	Self-adjointness implies that $\D_P-i$ is invertible.
	For $u=(\D_P-i)^{-1}f$, we have
	\[
	\|f\|_{L^2}^2
	=
	\|(\D_P-i)u\|_{L^2}^2
	=
	\|\D_Pu\|_{L^2}^2+\|u\|_{L^2}^2,
	\]
	since $\langle\D_Pu,u\rangle_{L^2}$ is real.
	Together with \eqref{regularity spectral}, this shows
	that $(\D_P-i)^{-1}:L^2\to H^1$ is bounded.
	The compact embedding $H^1\hookrightarrow L^2$
	therefore makes $(\D_P-i)^{-1}$ compact on $L^2$.
	
	Consequently, $\D_P$ has compact resolvent and admits
	an $L^2$-orthonormal eigenbasis
	$(\varphi_j)_{j\in\IN}$ with real eigenvalues
	$(\lambda_j)_{j\in\IN}$, counted with multiplicity.
	Each eigenvalue has finite multiplicity, and
	$|\lambda_j|\to\infty$.
	
	If $\D_P$ is invertible, its spectral gap is
	\[
	\mu
	\coloneq
	\operatorname{dist}(0,\spec(\D_P))
	=
	\min_j|\lambda_j|>0.
	\]
	The spectral theorem then gives
	\[
	\|\D_P^{-1}\|_{\mathcal L(L^2,L^2)}
	=
	\mu^{-1}.
	\]
	Combining this with \eqref{regularity spectral} yields
	a bounded inverse
	\[
	\D_P^{-1}:L^2(\Sigma\Omega)\longrightarrow
	H^1(\Sigma\Omega).
	\]
	The estimates for its $L^2\to H^1$ and $L^2\to L^q$
	operator norms are given in
	Lemma~\ref{lem:inverse-Lq}.

	\begin{remark}[Properties of the Dirac operator with
		boundary condition]
		For smooth boundaries, the chiral projections above
		define local elliptic self-adjoint realizations of
		the Dirac operator satisfying the regularity estimate
		\eqref{regularity spectral}; see
		\cite{Farinelli98,Raulot08,Baer12,Grosse-Uribe-vanDenBosch}.
		Invertibility is a separate assumption.
	\end{remark}

\section{The Main Theorem}
	
Let $(\Omega,g,\sigma)$ be a smooth compact Riemannian
spin manifold of dimension $n\geq2$ with smooth boundary
and spinor bundle $\Sigma\Omega$.
Let $P$ be a bounded linear boundary operator on
$H^1(\Sigma\Omega)$ such that the realization
\[
\D_P v=\D v,
\qquad
\dom(\D_P)
=
\{v\in H^1(\Sigma\Omega):Pv=0\},
\]
is self-adjoint, invertible, and regular. Thus, for
some constant $c_1>0$,
\[
\|v\|_{H^1}^2
\leq
c_1\bigl(\|v\|_{L^2}^2+\|\D_Pv\|_{L^2}^2\bigr),
\qquad v\in\dom(\D_P).
\]
Write $
\mu=\operatorname{dist}(0,\spec(\D_P))>0$
for its spectral gap.

Throughout this section, set
\[
p=\frac{2n}{n-1},
\qquad
q=2(p-1)=\frac{2(n+1)}{n-1},
\qquad
F(v)=|v|^{p-2}v.
\]
These exponents give
\[
F:L^q(\Sigma\Omega)\longrightarrow L^2(\Sigma\Omega),
\qquad
\|F(v)\|_{L^2}=\|v\|_{L^q}^{p-1}.
\]
Since $\D_P^{-1}$ maps $L^2$ boundedly into $H^1$
and $H^1\hookrightarrow L^q$, the composition
$\D_P^{-1}F$ is well defined on $L^q$.

For the quantitative inverse estimates, we use the
Gagliardo--Nirenberg interpolation exponent $
s=\frac n{n+1}, $
which satisfies
\begin{equation}
	\frac1q=\frac12-\frac{s}{n},
	\qquad
	\frac{s}{n}=\frac1{n+1}.
	\label{eq:critical-exponent}
\end{equation}

We perform the contraction in $L^q$, after writing
$u=w+\psi$, where $\psi$ is the harmonic representative
of the prescribed boundary datum.
The \say{correction} map
\[
w\longmapsto
\lambda\D_P^{-1}F(w+\psi)
\]
takes values in $\dom(\D_P)$.
Consequently, its fixed point has $H^1$ regularity
and satisfies the homogeneous boundary condition,
so that $Pu=P\psi$.
No trace operator on the contraction space $L^q$
is required.

\subsection{Auxiliary estimates}

	Set
	\[
	A_q=\|\D_P^{-1}\|_{\mathcal L(L^2,L^q)},
	\qquad
	B_1=\|\D_P^{-1}\|_{\mathcal L(L^2,H^1)}.
	\]
	
	\begin{lemma}[Inverse estimates in ordinary function spaces]
		\label{lem:inverse-Lq}
		The constants $A_q$ and $B_1$ are finite, and
		\begin{align*}
			\sqrt{\mu^{-2}+n^{-1}} &\leq B_1
			\leq\sqrt{c_1}\sqrt{1+\mu^{-2}},\\
			\frac{\vol_g(\Omega)^{-1/(n+1)}}{\mu} &\leq A_q
			\leq C_{\mathrm{GN}}\mu^{s-1}B_1^s
			\leq C_{\mathrm{GN}}c_1^{s/2}
			\frac{(1+\mu^2)^{s/2}}{\mu}.
		\end{align*}
	\end{lemma}
	
	\begin{proof}
		The spectral theorem gives
		\[
		\|\D_P^{-1}f\|_{L^2}\leq\mu^{-1}\|f\|_{L^2}.
		\]
		Applying the regularity estimate to $\D_P^{-1}f$ therefore yields
		\[
		\|\D_P^{-1}f\|_{H^1}^2
		\leq c_1(1+\mu^{-2})\|f\|_{L^2}^2.
		\]
		The Gagliardo--Nirenberg inequality gives
		\[
		\|\D_P^{-1}f\|_{L^q}
		\leq C_{\mathrm{GN}}
		\|\D_P^{-1}f\|_{L^2}^{1-s}\|\D_P^{-1}f\|_{H^1}^s
		\leq C_{\mathrm{GN}}\mu^{s-1}B_1^s\|f\|_{L^2},
		\]
		and therefore the upper bound via
		\[
		\mu^{s-1}(1+\mu^{-2})^{s/2}
		=
		\frac{(1+\mu^2)^{s/2}}{\mu}.
		\]
		
		For the lower bounds, choose an eigenspinor
		\[
		\D_P\varphi=\nu\varphi,\qquad |\nu|=\mu,\qquad \|\varphi\|_{L^2}=1.
		\]
		Then
		\[
		A_q\geq \|\D_P^{-1}\varphi\|_{L^q}
		=\frac{\|\varphi\|_{L^q}}{\mu}.
		\]
		Since \(q>2\), Hölder's inequality gives
		\[
		1=\|\varphi\|_{L^2}
		\leq \vol_g(\Omega)^{1/2-1/q}\|\varphi\|_{L^q}.
		\]
		Consequently,
		\[
		A_q\geq \frac{\vol_g(\Omega)^{1/q-1/2}}{\mu}
		=\frac{\vol_g(\Omega)^{-1/(n+1)}}{\mu}.
		\]
		Similarly,
		\[
		B_1\geq
		\|\D_P^{-1}\varphi\|_{H^1}
		=
		\frac{\|\varphi\|_{H^1}}{\mu}.
		\]
		The pointwise estimate
		$|\D\varphi|\leq\sqrt n\,|\nabla^S\varphi|$
		and the eigenvalue equation imply
		\[
		\mu^2
		=
		\|\D\varphi\|_{L^2}^2
		\leq n\|\nabla^S\varphi\|_{L^2}^2.
		\]
		Hence
		\[
		\|\varphi\|_{H^1}^2
		=
		1+\|\nabla^S\varphi\|_{L^2}^2
		\geq 1+\frac{\mu^2}{n},
		\]
		and therefore
		\[
		B_1\geq\sqrt{\mu^{-2}+\frac1n}.
		\]
	\end{proof}
			
	\begin{lemma}\label{lem:F}
		Let $F(v)=|v|^{p-2}v$ and $q=2(p-1)$.
		Then $F:L^q(\Sigma\Omega)\to L^2(\Sigma\Omega)$ is locally Lipschitz:
		\begin{equation}
			\|F(v)\|_{L^2}=\|v\|_{L^q}^{p-1},
			\label{F growth}
		\end{equation}
		and
		\begin{equation}
			\|F(v)-F(z)\|_{L^2}
			\leq
			(p-1)\max\{\|v\|_{L^q},\|z\|_{L^q}\}^{p-2}
			\|v-z\|_{L^q}.
			\label{F Lipschitz}
		\end{equation}
	\end{lemma}
	
	\begin{proof}
		The growth identity follows from $q=2(p-1)$.
		Regard each fiber as a real inner product space. For $\zeta\ne0$,
		\[
		DF(\zeta)h
		=
		|\zeta|^{p-2}h
		+
		(p-2)|\zeta|^{p-4}
		\operatorname{Re}\langle\zeta,h\rangle\zeta,
		\]
		and $DF(0)=0$. Hence
		\[
		|DF(\zeta)h|\leq(p-1)|\zeta|^{p-2}|h|.
		\]
		The fundamental theorem of calculus and H\"older's inequality, with
		\[
		\frac12=\frac{p-2}{q}+\frac1q,
		\]
		give
		\[
		\|F(v)-F(z)\|_{L^2}
		\leq
		(p-1)\int_0^1
		\|(1-t)z+tv\|_{L^q}^{p-2}\|v-z\|_{L^q}\,dt.
		\]
		The triangle inequality bounds the intermediate norm by
		$\max\{\|v\|_{L^q},\|z\|_{L^q}\}$, proving the result.
	\end{proof}

Since the boundary value problem only sees the values of $P\psi_0$, we are free to replace $\psi_0$ by a harmonic spinor $\Pi\psi_0$ with the same boundary value.

\begin{lemma}[Harmonic extension]\label{lem:harmonic-extension}
	For every $\psi_0\in H^1(\Sigma\Omega)$ there is a unique $\psi\in H^1(\Sigma\Omega)$ with $\D\psi=0$ and $P\psi=P\psi_0$, namely
	\[
	\psi=\psi_0-\D_P^{-1}(\D\psi_0).
	\]
	Moreover,
		\begin{align*}
			\|\psi\|_{L^q}
			&\leq\|\psi_0\|_{L^q}+A_q\|\D\psi_0\|_{L^2},\\
			\|\psi\|_{H^1}
			&\leq\|\psi_0\|_{H^1}+B_1\|\D\psi_0\|_{L^2}.
		\end{align*}
	
\end{lemma}

\begin{proof}
	$w\coloneq\D_P^{-1}(\D\psi_0)\in\dom(\D_P)=\ker(P)$, so $Pw=0$ and $P\psi=P\psi_0$. On $\dom(\D_P)$ the realization equals $\D$, so $\D w=\D\psi_0$, and hence $\D\psi=0$.
	If $\psi_1,\psi_2$ have the stated properties, their difference belongs to $\ker(P)$ and satisfies $\D(\psi_1-\psi_2)=0$, hence lies in $\ker\D_P=\{0\}$.
	Finally, $\psi\in H^1$ by regularity of the inverse.
	The estimates follow from the definitions of $A_q$ and $B_1$.
\end{proof}

\begin{remark}[Continuity of the harmonic extension]\label{rem:harmonic-continuity}
	The map
	\[
	\Pi:H^1(\Sigma\Omega)\to H^1(\Sigma\Omega),
	\qquad
	\Pi\psi_0=\psi_0-\D_P^{-1}(\D\psi_0),
	\]
	is a bounded linear projection onto the harmonic sections and preserves the value of $P$.
	Its value depends only on $P\psi_0$. By the pointwise bound $|\D\psi_0|\leq \sqrt{n} |\nabla^S\psi_0|$,
		\[
		\|\D\psi_0\|_{L^2}\leq\sqrt n\,\|\nabla^S\psi_0\|_{L^2},
		\]
		and hence
		\[
		\|\Pi\psi_0\|_{H^1}
		\leq(1+\sqrt n\,B_1)\|\psi_0\|_{H^1}.
		\]
\end{remark}

	From now on, let $\psi$ be harmonic and write $u=w+\psi$.
	The boundary value problem is equivalent to
	\[
	\D_Pw=\lambda F(w+\psi).
	\]
	Define
	\[
	T_\psi:L^q(\Sigma\Omega)\to L^q(\Sigma\Omega),
	\qquad
	T_\psi w=\lambda \D_P^{-1}F(w+\psi).
	\]
	For every $w\in L^q$,
	\[
	T_\psi w\in\dom(\D_P)\subset H^1.
	\]

\subsection{Existence and uniqueness}\label{thm:existence-uniqueness}

	\begin{theorem}[Existence and uniqueness via fixed point]
		\label{thm:fixed-point}
		Let $\psi\in H^1(\Sigma\Omega)$ be harmonic and set $m=\|\psi\|_{L^q}$.
		For $R>0$, define
		\[
		\mathcal B_R=\{w\in L^q(\Sigma\Omega):\|w\|_{L^q}\leq R\}.
		\]
		Assume
		\begin{align}
			|\lambda|A_q(R+m)^{p-1}&\leq R,\\
			\Theta\coloneq(p-1)|\lambda|A_q(R+m)^{p-2}&<1.
			\label{condition: contraction}
		\end{align}
		Then $T_\psi$ is a $\Theta$-contraction on $\mathcal B_R$.
		Its unique fixed point $w$ belongs to $\dom(\D_P)$, and $u=w+\psi$ solves
		\[
		\D u=\lambda|u|^{p-2}u,\qquad Pu=P\psi.
		\]
		The solution is unique among all $H^1$-solutions satisfying
		$\|u-\psi\|_{L^q}\leq R$.
		
		For every $w_0\in\mathcal B_R$, the iteration
		\[
		w_{k+1}=T_\psi w_k
		\]
		converges to $w$ in $L^q$. Every iterate with $k\geq1$ belongs to $\dom(\D_P)$.
	\end{theorem}
	
	\begin{proof}
		For $w\in\mathcal B_R$, the triangle inequality gives
		$\|w+\psi\|_{L^q}\leq R+m$. 
		From Lemma \ref{lem:F}, we know
		\[ \|F(w + \psi)\|_{L^2} = \|w + \psi\|_{L^q}^{p-1} \leq (R + m)^{p-1}.\]
		
		By definition, $A_q$ is the norm of the operator $\D_P^{-1}$ from $L^2$ to $L^q$, meaning we get
		\[ \|T_\psi w\|_{L^q} = \|\lambda \D_P^{-1} F(w + \psi)\|_{L^q} \leq |\lambda|  A_q  \|F(w + \psi)\|_{L^2}, \]
		and consequently,
		\[
		\|T_\psi w\|_{L^q}
		\leq|\lambda|A_q(R+m)^{p-1}\leq R.
		\]
		Thus $T_\psi$ is a self-map on the ball $\cB_R$.
		
		For $w,z\in\mathcal B_R$, we similarly have
		\[
		\|T_\psi w - T_\psi z\|_{L^q} \leq |\lambda| A_q \|F(w + \psi) - F(z + \psi)\|_{L^2}.
		\]
	 Using Lemma~\ref{lem:F}, this gives
		\[
		\|T_\psi w-T_\psi z\|_{L^q}
		\leq(p-1)|\lambda|A_q(R+m)^{p-2}\|w-z\|_{L^q}
		=
		\Theta\|w-z\|_{L^q}.
		\]
		The Banach Fixed-Point Theorem applies by the assumption $\Theta<1$ and  since $\mathcal B_R$ is complete, as it is a closed subset of the Banach-space $L^q(\Sigma \Omega)$.
		
		The fixed-point identity
		\[
		w=\lambda \D_P^{-1}F(w+\psi)
		\]
		implies $w\in\dom(\D_P)$ and $\D_Pw=\lambda F(w+\psi)$.
		Since $\D\psi=0$ and $Pw=0$, this proves the equation and boundary condition for $u=w+\psi$.
		
		Conversely, every $H^1$-solution in the stated ball yields a fixed point $w=u-\psi$. This proves the uniqueness assertion. The statement about the iteration is immediate by the Banach Fixed-Point Theorem and since $\t{Im}(T_\psi)\subset\dom(\D_P)$.
	\end{proof}

Choosing the radius $R$ as a function of $\psi$, we retrieve the following.

	\begin{corollary}[Adapted radius]\label{cor:adapted-radius}
		Assume $\psi\not=0$. For $\alpha>0$, define
		\[
		\delta_p(\alpha)
		=
		\min\left\{
		\frac{\alpha}{(1+\alpha)^{p-1}},
		\frac1{(p-1)(1+\alpha)^{p-2}}
		\right\}.
		\]
		Let $\psi\ne0$ be harmonic. If
		\begin{equation}
			|\lambda|A_q\|\psi\|_{L^q}^{p-2}<\delta_p(\alpha),
			\label{eq:smallness}
		\end{equation}
		then the boundary value problem
		\[
		\D u=\lambda|u|^{p-2}u \quad\text{in }\Omega,
		\qquad
		Pu=P\psi \quad\text{on }\partial\Omega
		\]
		has a solution $u\in H^1(\Sigma\Omega)$ satisfying
		\[
		\|u-\psi\|_{L^q}\leq\alpha\|\psi\|_{L^q}.
		\]
		It is unique among all $H^1$-solutions satisfying this bound.
		The iteration converges for every initial value in the corresponding ball, including $w_0=0$.
	\end{corollary}
	
	\begin{proof}
		Apply Theorem~\ref{thm:fixed-point} with $R=\alpha\|\psi\|_{L^q}\not=0$.
		The two terms defining $\delta_p(\alpha)$ give the self-map and contraction requirements.
	\end{proof}

\begin{remark}[The parameter $\alpha$ and the size of the solution]
	\label{rem:factor-two}
	
		Taking small $\alpha$ restricts the search to a small ball around the harmonic representative. This improves the Lipschitz estimate but makes the self-map condition more restrictive.
		Taking large $\alpha$ enlarges the ball while increasing the Lipschitz constant.
		
		The two terms in $\delta_p(\alpha)$ agree at
		\[
		\alpha_*=\frac1{p-2}.
		\]
		The first term is increasing on $(0,\alpha_*)$ and the second is decreasing on $(0,\i)$.
		Consequently,
		\[
		\max_{\alpha>0}\delta_p(\alpha)
		=
		\delta_p(\alpha_*)
		=\frac{(p-2)^{p-2}}{(p-1)^{p-1}}.
		\]
		For the spinorial Yamabe exponent, we have $\alpha_*=\frac{n-1}{2}.$
		This is the best sufficient threshold obtained by optimizing these two estimates. If one wants uniqueness of solutions on larger balls, or an existence within a smaller set, then the parameter $\alpha$ may be chosen differently to suit these goals.
	
\end{remark}

	Theorem~\ref{thm:main_stability} follows from Corollary~\ref{cor:adapted-radius} with $\alpha=\alpha_*$ and the convergence estimates below.
	
	\begin{corollary}[Sufficient conditions]
		\label{cor:spectral-condition}
		The conclusion of Corollary~\ref{cor:adapted-radius} holds if
		\[
		|\lambda|\|\psi\|_{L^q}^{p-2}
		<
		\frac{\delta_p(\alpha)}
		{C_{\mathrm{GN}}c_1^{s/2}}
		\frac{\mu}{(1+\mu^2)^{s/2}}.
		\]
	\end{corollary}
	
	\begin{proof}
		Use Lemma~\ref{lem:inverse-Lq} in \eqref{eq:smallness}.
	\end{proof}

\subsection{Quantitative convergence and $H^1$-control}

	\begin{theorem}[Quantitative convergence and $H^1$-control]
		\label{prop:quantitative}
		Assume the conditions of Theorem~\ref{thm:fixed-point}, and put
		\[
		L_R=(p-1)(R+m)^{p-2},
		\qquad
		\Theta=|\lambda|A_qL_R<1.
		\]
		Let $w$ be the fixed point. Then
		\[
		\|w_k-w\|_{L^q}
		\leq\frac{\Theta^k}{1-\Theta}\|w_1-w_0\|_{L^q},
		\qquad k\geq1.
		\]
		For $k\geq0$,
		\begin{align}
			\|\D_P(w_{k+1}-w)\|_{L^2}
			&\leq|\lambda|L_R\|w_k-w\|_{L^q},\\
			\|w_{k+1}-w\|_{H^1}
			&\leq|\lambda|B_1L_R\|w_k-w\|_{L^q},
			\label{estimate:H^1_convergence}
		\end{align}
		Thus the iteration converges geometrically in these spaces.
		
		Moreover,
		\[
		\|w_{k+1}\|_{H^1},\ \|w\|_{H^1}
		\leq|\lambda|B_1(R+m)^{p-1}.
		\]
		For $u_k=w_k+\psi$ and $k\geq1$, this gives
		\[
		\|u_k\|_{H^1}
		\leq\|\psi\|_{H^1}+|\lambda|B_1(R+m)^{p-1}.
		\]
	\end{theorem}
	
	\begin{proof}
		The $L^q$ estimate is the standard contraction error bound.
		For the remaining estimates, use
		\[
		w_{k+1}-w
		=
		\lambda \D_P^{-1}\bigl(F(w_k+\psi)-F(w+\psi)\bigr)
		\]
		and
		\[
		\|F(w_k+\psi)-F(w+\psi)\|_{L^2}
		\leq L_R\|w_k-w\|_{L^q}.
		\]
		Apply $\D_P$ and the boundedness of $\D_P^{-1}:L^2\to H^1$ respectively to get the two estimates.
		The uniform upper bounds then follow from
		$\|F(w_k+\psi)\|_{L^2}\leq(R+m)^{p-1}$ and the definition of $B_1$, as well as triangle inequality.
	\end{proof}

\subsection{Continuous dependence on the boundary data}

	\begin{theorem}[Continuous dependence on the boundary data]
		\label{prop:stability}
		Let $\alpha>0$ and $M>0$ satisfy
		\[
		|\lambda|A_qM^{p-2}<\delta_p(\alpha).
		\]
		Set
		\[
		R=\alpha M,\qquad
		L_M=(p-1)(1+\alpha)^{p-2}M^{p-2},
		\qquad
		\Theta=|\lambda|A_qL_M<1.
		\]
		Let $\psi_1,\psi_2\in H^1(\Sigma\Omega)$ be harmonic with
		$\|\psi_j\|_{L^q}\leq M$.
		For $j=1,2$, let $w_j$ be the unique fixed point of
		\[
		T_{\psi_j}(w)=\lambda\D_P^{-1}F(w+\psi_j)
		\]
		in the common ball
		\[
		\mathcal B_R
		=
		\{w\in L^q(\Sigma\Omega):\|w\|_{L^q}\leq R\},
		\]
		and put $u_j=w_j+\psi_j$. Then
		\begin{align*}
			\|w_1-w_2\|_{L^q}
			&\leq\frac{\Theta}{1-\Theta}\|\psi_1-\psi_2\|_{L^q},\\
			\|u_1-u_2\|_{L^q}
			&\leq\frac1{1-\Theta}\|\psi_1-\psi_2\|_{L^q},\\
			\|w_1-w_2\|_{H^1}
			&\leq\frac{|\lambda|B_1L_M}{1-\Theta}
			\|\psi_1-\psi_2\|_{L^q}.
		\end{align*}
		In particular,
		\[
		\|u_1-u_2\|_{H^1}
		\leq
		\|\psi_1-\psi_2\|_{H^1}
		+
		\frac{|\lambda|B_1L_M}{1-\Theta}
		\|\psi_1-\psi_2\|_{L^q}.
		\]
	\end{theorem}
	
	\begin{proof}
		For $j=1,2$ and $w\in\mathcal B_R$, we have
		\[
		\|w+\psi_j\|_{L^q}
		\leq R+M=(1+\alpha)M.
		\]
		The smallness assumption therefore gives
		\[
		|\lambda|A_q(R+M)^{p-1}<R,
		\qquad
		(p-1)|\lambda|A_q(R+M)^{p-2}=\Theta<1.
		\]
		Thus Theorem~\ref{thm:fixed-point} applies to both maps
		$T_{\psi_j}$ on the same ball $\mathcal B_R$.
		Their fixed-point identities yield
		\[
		w_1-w_2=\lambda \D_P^{-1}(F(u_1)-F(u_2)).
		\]
		Since $\|u_j\|_{L^q}\leq(1+\alpha)M$,
		\[
		\|w_1-w_2\|_{L^q}
		\leq\Theta\|u_1-u_2\|_{L^q}.
		\]
		Combine this with the triangle inequality
		\[
		\|u_1-u_2\|_{L^q}
		\leq\|w_1-w_2\|_{L^q}+\|\psi_1-\psi_2\|_{L^q}
		\]
		and absorb the term with coefficient $\Theta$ to get  $\|w_1-w_2\|_{L^q}
		\leq\frac{\Theta}{1-\Theta}\|\psi_1-\psi_2\|_{L^q}$. Inserting this into our $\|u_1-u_2\|_{L^q}$ estimate yields the second inequality $\|u_1-u_2\|_{L^q}
		\leq\frac1{1-\Theta}\|\psi_1-\psi_2\|_{L^q}$.
		The remaining bounds follow by
		\begin{align*} \|w_1 - w_2\|_{H^1} &= \|\lambda \D_P^{-1}(F(u_1) - F(u_2))\|_{H^1} \\
	 &\leq |\lambda| B_1 \|F(u_1) - F(u_2)\|_{L^2}  \\
	 &\leq |\lambda| B_1 L_M \|u_1 - u_2\|_{L^q} 
		\end{align*}
		and triangle inequality.
	\end{proof}

\begin{remark}
	Theorem~\ref{thm:fixed-point} provides an iteration that may serve as the basis of a numerical method, as it yields existence and uniqueness in a specified ball, along with a sequence constructed by solving a linear boundary value problem in each iteration step.
	Theorems~\ref{prop:quantitative} and~\ref{prop:stability} show that this iteration is stable and converges geometrically.
	For $k\geq1$, the a posteriori estimate
		\[
		\|w_k-w\|_{L^q}
		\leq\frac{\Theta}{1-\Theta}\|w_k-w_{k-1}\|_{L^q}
		\]
		controls the error through successive iterates. Together with \eqref{estimate:H^1_convergence}, it also controls the $H^1$ error.
		
		Note that convergence of a discretized scheme additionally requires stability and consistency estimates for the discretization and control of linear-solver errors.
\end{remark}

\subsection{Further remarks on the method}

\begin{remark}[Homogeneity and scaling]
	
		The transformation
		\[
		(u,\psi,\lambda)\mapsto(tu,t\psi,t^{2-p}\lambda),
		\qquad t>0,
		\]
		preserves the equation and the form of its boundary data.
		The quantity $|\lambda|A_q\|\psi\|_{L^q}^{p-2}$ is invariant.
		
		Likewise,
		\[
		(\D_P,\lambda)\mapsto(a\D_P,a\lambda),
		\qquad a>0,
		\]
		leaves the boundary value problem unchanged. Since
		\[
		A_q(a\D_P)=a^{-1}A_q(\D_P),
		\]
		the smallness condition is invariant under this scaling too.
		
		The estimate
		\[
		A_q\leq C_{\mathrm{GN}}\mu^{s-1}B_1^s
		\]
		has the same homogeneity: $\mu$ scales by $a$ and $B_1$ by $a^{-1}$.
	
\end{remark}

\begin{remark}\label{remark: generality of method}
	The spinorial Yamabe equation was the motivation, but the fixed-point proof extends in several directions.
	\begin{enumerate}
		\item \emph{Other operators.}
		The argument applies to other operators on Hermitian bundles whenever the homogeneous boundary problem has a bounded inverse $\D_P^{-1}:L^2\to H^1$. Self-adjointness is not required for the fixed-point,
		convergence, or stability arguments once the bounded
		inverse $L_P^{-1}:L^2\to H^1$ is available.
		In the self-adjoint setting, it additionally provides
		the spectral characterization of the $L^2$ inverse
		norm used in our quantitative estimates.
		For non-self-adjoint realizations, the required inverse
		bounds must be established separately.
		
		\item \emph{Other exponents.}
		The proof applies for every $p>2$ such that
			\[
			H^1(E)\hookrightarrow L^{2(p-1)}(E).
			\]
			This includes
			\[
			2<p\leq\frac{2n-2}{n-2}\quad(n>2),
			\]
			and every finite $p>2$ when $n=2$.
			Set
			\[
			q=2(p-1),\qquad
			s_p=\frac{n(p-2)}{2(p-1)}.
			\]
			Then $1/q=1/2-s_p/n$.
			
			The restriction $p>2$ comes from the nonlinear
			contraction argument, rather than from the Sobolev
			embedding alone. For $1<p<2$, the map
			$F(v)=|v|^{p-2}v$, with $F(0)=0$, is H\"older
			continuous of exponent $p-1$, but is not locally
			Lipschitz at the origin. Indeed, for $v\neq0$ and $t>0$,
			\[
			\frac{|F(tv)-F(0)|}{|tv|}
			=
			t^{p-2}|v|^{p-2}\longrightarrow\infty
			\qquad(t\downarrow0).
			\]
			Thus the Lipschitz estimate needed for our contraction
			argument is unavailable. Moreover, if $1<p<3/2$, then
			$q=2(p-1)<1$, so the chosen $L^q$ space is not a
			Banach space with its usual expression for the norm.
			These observations do not rule out existence for
			$1<p<2$, but this range requires special arguments.
			
			For $n>2$, the endpoint
			\[
			p=\frac{2n-2}{n-2},
			\qquad
			q=\frac{2n}{n-2},
			\qquad
			s_p=1
			\]
			is included in the fixed-point argument, since the
			continuous Sobolev embedding
			$H^1(E)\hookrightarrow L^{2n/(n-2)}(E)$ remains valid.
			At this endpoint, one estimates
			\[
			A_q\leq C_{\mathrm S}B_1
			\]
			directly by the Sobolev inequality, rather than invoking
			Lemma~\ref{Fractional Gagliardo--Nirenberg} with $s<1$.
			
			The elementary higher-integrability bootstrap used in
			Section~\ref{Regularity via Bootstrapping} does not
			yield an improvement at this endpoint.
			Indeed, starting from $u\in H^1$, it gives only
			\[
			F(u)\in
			L^{\frac{2n}{(n-2)(p-1)}}=L^2,
			\]
			and hence returns to $H^1$.
			Higher regularity at the endpoint therefore requires
			an additional argument.
			
			For $p=2$, the problem is linear and no fixed-point
			argument is needed. Writing $u=w+\psi$, with $\psi$
			the harmonic representative of the boundary datum,
			gives
			\[
			(\D_P-\lambda)w=\lambda\psi.
			\]
			Thus, whenever $\lambda\notin\operatorname{spec}(\D_P)$,
			the unique solution is
			\[
			u=\psi+\lambda(\D_P-\lambda)^{-1}\psi,
			\]
			without any smallness condition on the boundary datum.
			
		\item \emph{Variable coefficients.}
		Let $k\in L^\infty(\Omega)$ and consider
		\[
		\D u=k(x)|u|^{p-2}u,\qquad Pu=P\psi.
		\]
		Replacing $|\lambda|$ by $\|k\|_{L^\infty}$ in the growth and Lipschitz estimates gives the same fixed-point conclusions, i.e. existence, uniqueness, convergence, and stability.
		\item \emph{Other non-linearities.}
		Let $N:L^q(E)\to L^2(E)$ be locally Lipschitz and
		consider the boundary value problem
		\[
		\D u=N(u)\quad\text{in }\Omega,
		\qquad
		Pu=P\psi_0\quad\text{on }\partial\Omega.
		\]
		Let $\psi$ be the harmonic representative of the
		boundary datum, so that
		\[
		\D\psi=0,\qquad P\psi=P\psi_0.
		\]
		Writing $u=w+\psi$, the problem is equivalent to
		the fixed-point equation
		\[
		w=T_\psi(w),
		\qquad
		T_\psi(w)=\D_P^{-1}N(w+\psi).
		\]
		
		Suppose that, for some $R>0$, there are constants
		$M_R,L_R\geq0$ such that
		\[
		\|N(w+\psi)\|_{L^2}\leq M_R
		\]
		and
		\[
		\|N(w+\psi)-N(z+\psi)\|_{L^2}
		\leq L_R\|w-z\|_{L^q}
		\]
		whenever $\|w\|_{L^q},\|z\|_{L^q}\leq R$.
		Then
		\[
		A_qM_R\leq R,\qquad A_qL_R<1
		\]
		ensure that $T_\psi$ maps the closed ball
		$\mathcal B_R\subset L^q(E)$ into itself and is
		a contraction there.
		Consequently, the boundary value problem has a
		unique $H^1$-solution among those satisfying
		$\|u-\psi\|_{L^q}\leq R$.
		The iteration converges geometrically in $L^q$
		and, after the first step, in $H^1$.
		Continuous dependence on the harmonic representative
		follows under uniform growth and Lipschitz estimates
		on the relevant translated balls.
		Higher regularity requires additional assumptions on $N$.
		
		As a demonstration, an additive source term is included in this framework
		by taking
		\[
		N(v)=\lambda|v|^{p-2}v+f,
		\qquad f\in L^2(E).
		\]
		The boundary value problem then reads
		\[
		\D u=\lambda|u|^{p-2}u+f,
		\qquad Pu=P\psi_0,
		\]
		and, with $q=2(p-1)$ and $m=\|\psi\|_{L^q}$,
		one may take
		\[
		M_R=|\lambda|(R+m)^{p-1}+\|f\|_{L^2},
		\qquad
		L_R=(p-1)|\lambda|(R+m)^{p-2}.
		\]
		Notice that the self-mapping condition gets more restrictive, but the condition: contraction remains invariant here.
		Estimating its contribution after applying the inverse
		gives the sharper sufficient self-map condition
		\[
		|\lambda|A_q(R+m)^{p-1}
		+\|\D_P^{-1}f\|_{L^q}\leq R.
		\]
		The condition: contraction remains
		\[
		(p-1)|\lambda|A_q(R+m)^{p-2}<1,
		\]
		because $f$ cancels when taking differences.
		
		In Section~\ref{section:brezis--Nirenberg}, we will write out the details for the Brezis--Nirenberg problem, associated to the function $N(v)=\lambda |v|^{p-2}v + av$.
	\end{enumerate}
\end{remark}

\begin{remark}[Auxiliary shifts and the non-invertible case]
	\label{rem:a-f}
	Let $\D_P$ be self-adjoint and regular, but not
	necessarily invertible.
	For an extension $\psi_0\in H^1(\Sigma\Omega)$
	of the prescribed boundary datum and
	$a\in\mathbb C\setminus\operatorname{spec}(\D_P)$,
	define
	\[
	T_{\psi_0,a}(w)
	=
	(\D_P-a)^{-1}
	\bigl(
	\lambda F(w+\psi_0)-\D\psi_0-aw
	\bigr).
	\]
	The associated iteration is
	\[
	w_{k+1}
	=
	(\D_P-a)^{-1}
	\bigl(
	\lambda F(w_k+\psi_0)-\D\psi_0-aw_k
	\bigr),
	\qquad k\geq0.
	\]
	Starting from $w_0\in L^q$, every iterate with $k\geq1$
	belongs to $\dom(\D_P)$.
	
	Equivalently, writing $u_k=w_k+\psi_0$, one obtains
	\[
	u_{k+1}
	=
	\psi_0+
	(\D_P-a)^{-1}
	\bigl(
	\lambda F(u_k)-\D\psi_0-a(u_k-\psi_0)
	\bigr).
	\]
	Thus each step solves the linear boundary value problem
	\[
	(\D-a)u_{k+1}
	=
	\lambda F(u_k)-au_k,
	\qquad
	Pu_{k+1}=P\psi_0.
	\]
	The term $-au_k$ compensates for the auxiliary shift,
	so any fixed point solves the original, unshifted
	nonlinear equation.
	
	Since $L^q\hookrightarrow L^2$ and
	$F:L^q\to L^2$, this map $T_{\psi_0,a}$ takes $L^q$ into
	$\dom(\D_P)\subset H^1$.
	Its fixed points correspond exactly to solutions
	$u=w+\psi_0$ of
	\[
	\D u=\lambda F(u),
	\qquad Pu=P\psi_0.
	\]
	Indeed, the fixed-point equation gives
	\[
	(\D_P-a)w
	=
	\lambda F(w+\psi_0)-\D\psi_0-aw,
	\]
	and cancellation of the terms $-aw$ yields the
	original equation.
	
	Set
	\[
	C_a=\|(\D_P-a)^{-1}\|_{\mathcal L(L^2,L^q)},
	\qquad
	K_\Omega=\operatorname{vol}_g(\Omega)^{1/2-1/q},
	\qquad
	m_0=\|\psi_0\|_{L^q}.
	\]
	By H\"older's inequality,
	\[
	\|w\|_{L^2}\leq K_\Omega\|w\|_{L^q}.
	\]
	Consequently, on the ball
	$\mathcal B_R=\{w:\|w\|_{L^q}\leq R\}$,
	\[
	\begin{aligned}
		\|T_{\psi_0,a}(w)\|_{L^q}
		&\leq
		C_a\bigl(
		|\lambda|(R+m_0)^{p-1}
		+|a|K_\Omega R
		\bigr)\\
		&\quad+
		\|(\D_P-a)^{-1}\D\psi_0\|_{L^q}.
	\end{aligned}
	\]
	Moreover,
	\[
	\|T_{\psi_0,a}(w)-T_{\psi_0,a}(z)\|_{L^q}
	\leq
	\Theta_a\|w-z\|_{L^q},
	\]
	where
	\[
	\Theta_a
	=
	C_a\bigl(
	(p-1)|\lambda|(R+m_0)^{p-2}
	+|a|K_\Omega
	\bigr).
	\]
	Thus a sufficient condition for contraction on
	$\mathcal B_R$ is that the preceding self-map
	bound be at most $R$ and that $\Theta_a<1$.
	
	These estimates cannot establish contraction when
	$\ker(\D_P)\neq\{0\}$.
	Indeed, the spectral theorem and the embedding
	$L^q\hookrightarrow L^2$ give
	\[
	\frac{1}{\operatorname{dist}
		(a,\operatorname{spec}(\D_P))}
	=
	\|(\D_P-a)^{-1}\|_{\mathcal L(L^2,L^2)}
	\leq K_\Omega C_a.
	\]
	If $0\in\operatorname{spec}(\D_P)$, then
	$\operatorname{dist}(a,\operatorname{spec}(\D_P))
	\leq|a|$, so
	\[
	|a|K_\Omega C_a\geq1,
	\qquad
	\Theta_a\geq1.
	\]
	This shows that the above sufficient contraction
	criterion fails; it does not, by itself, rule out
	convergence of the shifted iteration.
	
	Furthermore, every solution must satisfy the
	compatibility conditions
	\[
	\langle
	\lambda F(u)-\D\psi_0,h
	\rangle_{L^2}=0
	\qquad
	\text{for every }h\in\ker(\D_P).
	\]
	An auxiliary shift leaves these conditions unchanged.
	A treatment of the kernel therefore requires
	additional arguments, for example a decomposition
	\[
	L^2(\Sigma\Omega)
	=
	\ker(\D_P)\oplus\ker(\D_P)^\perp
	\]
	and an analysis of the resulting equations.
	
	Unlike the linear perturbation in
	Corollary~\ref{cor:linear-perturbation}, the shift
	considered here does not change the differential
	equation: it introduces the compensating term $-aw$
	into the iteration.
\end{remark}
	
\begin{remark}[Choice of the boundary extension]
	\label{non-harmonic}
	Assume that $\D_P$ is invertible.
	For an arbitrary extension $\psi_0\in H^1$ of
	the prescribed boundary datum, the iteration for
	$u=w+\psi_0$ is
	\[
	T_{\psi_0}(w)
	=
	\D_P^{-1}
	\bigl(
	\lambda F(w+\psi_0)-\D\psi_0
	\bigr).
	\]
	Set
	\[
	m_0=\|\psi_0\|_{L^q},
	\qquad
	h=\D_P^{-1}\D\psi_0,
	\qquad
	\varepsilon=\|h\|_{L^q}.
	\]
	On the ball $\|w\|_{L^q}\leq R$, sufficient
	self-map and condition: contractions are
	\[
	|\lambda|A_q(R+m_0)^{p-1}+\varepsilon\leq R,
	\]
	and
	\[
	(p-1)|\lambda|A_q(R+m_0)^{p-2}<1.
	\]
	The additive term $\varepsilon$ must therefore
	be included when using this extension and a ball
	centered at zero in the correction variable.
	
	Passing to the harmonic representative
	$\psi=\psi_0-h$ removes the explicit forcing term.
	More precisely, the corresponding map
	\[
	T_\psi(z)=\lambda\D_P^{-1}F(z+\psi)
	\]
	satisfies
	\[
	T_\psi(w+h)=T_{\psi_0}(w)+h.
	\]
	Thus the two iterations are related by translation;
	with corresponding initial values, they generate
	the same sequence of sections $u_k$.
	Their contraction properties on corresponding
	translated sets are consequently equivalent.
	
	The harmonic representative provides a canonical
	choice depending only on the prescribed boundary
	datum and simplifies the sufficient estimates.
	No minimization of the $L^q$ norm among extensions,
	or optimality among all choices of contraction balls,
	is asserted.
\end{remark}

\begin{remark}\label{rem:not smooth}
	The boundary of $\Omega$ does not need to be smooth,
	as the contraction theorem relies principally on the relevant
	Sobolev embedding
	\[
	H^1(\Sigma\Omega)\hookrightarrow
	L^{2(p-1)}(\Sigma\Omega)
	\]
	and a bounded inverse $L^2\to H^1$ for the chosen
	realization.
	
	The Gagliardo--Nirenberg inequality used to estimate
	$A_q$ also holds on a large class of nonsmooth domains,
	including bounded Lipschitz domains, where it follows
	from suitable extension operators and localization.
	This inequality is not needed for the contraction
	argument itself once the boundedness of the inverse
	$L^2\to L^{2(p-1)}$ is known.
	
	On nonsmooth domains, the required embedding and
	operator regularity properties must be established
	separately.
	For suitable planar corner domains and boundary
	conditions, such questions are studied in
	\cite{VanDenBosch21}.
	For higher-dimensional questions on self-adjointness
	and regularity of Dirac operators with nonsmooth
	boundary, we refer to \cite{Pankrashkin25}.
	The embedding and regularity constants may depend
	on the corner geometry.
	This extension concerns existence and convergence
	in the stated spaces; it does not imply smoothness
	at corners.
\end{remark}

\begin{remark}
	If we take for example a square $\Omega_1$ in $\IR^2$, then
	\cite{VanDenBosch21} yields self-adjointness and regularity for the Dirac
	operator under (constant) generalized MIT boundary condition. As explained in
	Remark~\ref{rem:not smooth} our theorem applies to $\Omega_1$ even though it is not
	smooth: if the kernel of $\D_{P}^{\Omega_1}$ is trivial, then there is a solution
	to the spinorial Yamabe equation over $\Omega_1$ under chiral bag boundary condition
	for a certain smallness of some parameters. We can do that for an adjacent square
	$\Omega_2$ that shares an edge with $\Omega_1$, giving us two solutions that do not
	need to agree along the common edge. Whether a global solution can be synthesized
	from these local solutions remains an open question, potentially requiring domain
	decomposition techniques for first-order systems. Clearly, if the domain
	$\Omega_1\cup\Omega_2$ admits a self-adjoint, regular and invertible $\D_P$ (which
	rectangles do by \cite{VanDenBosch21}), then we can construct a solution via our
	main theorem under certain smallness conditions, but the combined rectangle
	$\Omega_1\cup\Omega_2$ may potentially have larger eigenvalue
	$|\la_1(\D_P)|$ than its constituents. Also note that we could have chosen the two
	squares to overlap and tried the same. 
\end{remark}

\begin{remark}[Homogeneous boundary data and closed manifolds]
	\label{rem:closed}
		It is important to note that the iterative scheme presented here is designed for small-data problems and cannot, in general, be used to construct nonzero solutions with homogeneous boundary data $P\psi_0=0$. Given $P\psi_0=0$, the homogenous representative is $\psi=0$ by invertibility of $\D_P$. Theorem~\ref{thm:fixed-point} then applies with $m=0$ and $R$ sufficiently small, but the unique fixed point (and hence limit of any iteration scheme in the contraction ball) is the trivial solution. Unfortunately, we cannot prove any non-trivial existence results in this case. Also, the numerical application in this case is questionable, because the contraction property breaks down for the homogenous problem:
		
		To see this, consider the map $T_0(v) = \lambda \D_P^{-1}F(v)$. If $u \neq 0$ is a solution to the homogeneous problem, then $T_0(u) = u$. However, due to the homogeneous boundary condition, we have $T_0(tu) = t^{p-1}u$ for all $t > 0$. Thus, for $k\in\IN$, \[T_0^k(tu)=t^{(p-1)^k}u.\]
		For $0<t<1$, this converges to zero, and for $t>1$, its norm tends to infinity. Thus, the homogeneous boundary problem cannot be treated with the same contraction argument and small deviations from the solutions diverge. A normalization at each step could potentially fix this behaviour in application. We can only formulate a norm-gap result:
	
\end{remark}

\begin{corollary} \label{cor:Lq-gap}
Let $u \in \dom(\D_P)$ be a non-trivial solution to the homogeneous problem $\D u = \lambda |u|^{p-2}u$ with $Pu=0$. Then $\|u\|_{L^q} \geq (|\lambda| A_q)^{-1/(p-2)}$, where $q=2(p-1)$ and $A_q = \|\D_P^{-1}\|_{\mathcal{L}(L^2, L^q)}$.
\end{corollary}

\begin{proof}
Since $u \in \dom(\D_P)$ and $\D_P$ is invertible, the equation $\D_P u = \lambda |u|^{p-2}u$ is equivalent to the fixed-point identity
\[ u = \lambda \D_P^{-1} \left( |u|^{p-2}u \right). \]
Taking the $L^q$-norm of both sides and applying the definition of the operator norm $A_q$, we obtain
\[ \|u\|_{L^q} = \| \lambda \D_P^{-1} ( |u|^{p-2}u ) \|_{L^q} \leq |\lambda| A_q \| |u|^{p-2}u \|_{L^2}. \]
By the definition of $q=2(p-1)$, the $L^2$-norm of the non-linearity satisfies
\[ \| |u|^{p-2}u \|_{L^2} = \left( \int_\Omega |u|^{2(p-1)} \, d\vol_g \right)^{1/2} = \|u\|_{L^q}^{p-1}. \]
Substituting this into the inequality yields
\[ \|u\|_{L^q} \leq |\lambda| A_q \|u\|_{L^q}^{p-1}. \]
Since $u \neq 0$, we divide by $\|u\|_{L^q}$ to find
\[ 1 \leq |\lambda| A_q \|u\|_{L^q}^{p-2}, \]
which rearranges to the desired lower bound
\[ \|u\|_{L^q} \geq \left( \frac{1}{|\lambda| A_q} \right)^{\frac{1}{p-2}}. \]
\end{proof}
\begin{remark}
	The statement is also true for closed manifolds, i.e. in the setting without any boundary (condition), as long as $\D$ has the necessary properties of being invertible and regular.
\end{remark}

\vspace{1cm}
\section{Regularity via Bootstrapping and A Priori Estimates}\label{Regularity via Bootstrapping}
Our goal is now, to employ a bootstrapping argument, in order to deduce higher regularity for solutions $u \in H_\cT^1(\Sigma \Omega)$ to the equation
\[
\D u = \lambda |u|^{p-2}u
\]
under the boundary condition $Pu = P\psi$. To successfully carry out this procedure, we require a suitable linear regularity statement together with an \emph{a priori} bound on $u$.

Such a regularity statement and \emph{a priori} bound always hold locally in the interior of the domain (see the upcoming Remark~\ref{interior a priori}). However, if we wish to achieve higher regularity on the whole compact set $\Omega$ including its boundary, we need them to be valid up to the boundary. In general, this is only possible for a suitable class of boundary operators $P$:

More precisely, we require that for every $q \in (1, \infty)$,
\[
\varphi\in H^1(\Sigma\Omega),\qquad
P\varphi=0,\qquad
\varphi,\D\varphi\in L^q(\Sigma\Omega)
\quad\Longrightarrow\quad
\varphi\in W^{1,q}(\Sigma\Omega),
\]
and that there exists a constant $c_q > 0$ such that the \emph{a priori} bound
\begin{equation}
	\| \varphi \|_{W ^{1,q}}^2 \leq c_{q} \Big( \|\varphi \|_{L^q}^2 + \| \D  \varphi \|_{L^q}^2 \Big)
	\label{L^p estimate}
\end{equation}
holds for these sections. Here $\D\varphi$ is understood distributionally and $P\varphi=0$ in the $H^1$ trace sense. On $W^{1,q}$ sections, the boundary operator extends boundedly as
\[
P\colon W^{1,q}_\cT(\Sigma\Omega) \to W^{1-1/q, q}_\cT\big(\restr{\Sigma \Omega}{\partial \Omega}\big),
\]
with values in the appropriate boundary subbundle when $P$ is a projection.

A standard sufficient condition for this regularity statement and the validity of~\eqref{L^p estimate}, for smooth local $P$ on a smooth domain, is that the boundary condition is elliptic in the sense of Shapiro--Lopatinski. Under this condition, the corresponding Sobolev regularity follows from elliptic boundary theory for all $q\in(1,\infty)$; see~\cite{Schwarz95}. Fredholmness alone should not be identified with the weak regularity implication stated above.

\setcounter{section}{1}
\setcounter{theorem}{1}
\begin{theorem}[Regularity]
	\label{second main theorem}
	Let $(\Omega,g)$ be a smooth compact Riemannian spin manifold of dimension $n\geq2$ with smooth boundary.
	Suppose that $u\in H^1(\Sigma\Omega)$ satisfies
	\[
	\D u=\lambda|u|^{p-2}u,
	\qquad
	p=\frac{2n}{n-1}.
	\]
	Then $u\in C_{\mathrm{loc}}^{1,\alpha}\bigl(\Sigma(\operatorname{int}\Omega)\bigr)$ for every $0<\alpha<1$, and $u$ is smooth in the interior away from its zero set.
	
	Furthermore, assume that $Pu=P\psi$ and that the boundary operator $P$ satisfies the \textbf{Shapiro--Lopatinski condition} for $\D$. Then:
	\begin{enumerate}
		\item If the prescribed boundary datum $\psi$ admits an extension in $\bigcap_{1<a<\infty} W^{1,a}(\Sigma\Omega)$, then $u\in \bigcap_{1<a<\infty} W^{1,a}(\Sigma\Omega)$.
		\item If $P$ additionally admits boundary Schauder estimates and $\psi$ admits a $C^{1,\alpha}$ extension, then $u\in C^{1,\alpha}(\Sigma\Omega)$.
		\item If the datum $\psi$ is smooth and $P$ admits boundary Schauder regularity of all orders, then $u$ is smooth away from its zero set up to the boundary.
	\end{enumerate}
\end{theorem}
\setcounter{section}{5}
\setcounter{theorem}{0}
\begin{remark}
	Here $\psi$ may be any extension of the prescribed boundary datum with the required regularity; it need not be harmonic. If we only assume that $\psi\in L^\i(\Sigma\Omega)$ with $\D\psi=0$, this does not imply $\psi\in W^{1,\i}$ up to the boundary. The homogeneous-boundary estimate must be applied to $u-\psi$, with the required regularity of $\psi$ assumed separately.
\end{remark}

\begin{remark}\label{interior a priori}
	If one does not have the regularity statement and \emph{a priori} bound for a specific boundary condition (i.e. not assuming the Shapiro--Lopatinski condition), we can still use their interior versions, which do not need any requirements on the boundary condition. For relatively compact open sets $U\Subset V\Subset\operatorname{int}(\Omega)$, interior elliptic regularity gives
	\begin{equation*}
		\| \varphi\|_{W ^{1,q}(\Sigma U)}^2
		\leq c(U,V,q) \left(
		\|\varphi \|_{L^q(\Sigma V)}^2
		+\| \D \varphi \|_{L^q(\Sigma V)}^2
		\right)
	\end{equation*}
	together with $\varphi\in W^{1,q}(\Sigma U)$ whenever $\varphi,\D\varphi\in L^q(\Sigma V)$; see the interior estimates in \cite[Thm.~3.2.1]{Ammann03hab}. These are local statements and require no global simple-connectedness assumption.
	
	Applying the bootstrap below on nested relatively compact sets, followed by interior Schauder regularity, shows that $u$ has a $C^{1,\alpha}_{\mathrm{loc}}$ representative for every $0<\alpha<1$, smooth in the interior away from $u^{-1}(\{0\})$. These interior conclusions require neither the boundary regularity assumptions nor a regular extension of the boundary datum.
\end{remark}
\begin{remark}\label{remark: Lq_ineq}
	For every $1<q<\infty$ and
	$\varphi\in W^{1,q}(\Sigma\Omega)$, the pointwise estimate
	\[
	|\D\varphi|\leq\sqrt n\,|\nabla^S\varphi|
	\]
	gives
	\[
	\|\varphi\|_{L^q}^2+\|\D\varphi\|_{L^q}^2
	\leq
	\|\varphi\|_{L^q}^2
	+n\|\nabla^S\varphi\|_{L^q}^2
	\leq
	(n+1)\|\varphi\|_{W^{1,q}}^2.
	\]
	No boundary condition is required.
\end{remark}

\begin{proof}[Proof of Thm. 1.2.:]
	First, assume $n>2$.
	We want to prove higher regularity of our solution $u\in H^1(\Sigma \Omega)$:
	
	Recall that the inclusion $L^s(\Sigma \Omega)\subset L^t(\Sigma \Omega)$ is bounded for all $s>t$ by the finite volume of $\Omega$ (compact) and H\"older's inequality $\Vert f \Vert_{L^t} \leq \operatorname{vol}_g(\Omega)^{\f1t - \f1s} \Vert f \Vert_{L^s} $. This further implies a bounded inclusion $W^{1,s}(\Sigma \Omega)\subset W^{1,t}(\Sigma \Omega)$ of Sobolev spaces for all $s>t$.
	
	Let $r\coloneq 2<n$, and notice that $u\in H^1(\Sigma\Omega) = W^{1,r}(\Sigma\Omega)$. We may apply the Sobolev embedding (Thm.~\ref{sobolev embeddings} (i)), to see that $u\in L^l(\Sigma \Omega)$, where $l$ is the solution of
	\begin{equation*}
		\frac{1}{l} = \frac{1}{r} - \frac{1}{n},
	\end{equation*}
	i.e. $l = \frac{rn}{n-r} = \frac{2n}{n-2}$.
	
	With $u\in L^l(\Sigma \Omega)$, the partial differential equation then implies $\D u \in L^{\frac{l}{p-1}}(\Sigma \Omega) = L^{r'}(\Sigma \Omega)$ with $r' \coloneq \frac{l}{p-1}$. We hence have $u\in L^{r'}$ by the inclusion property from above since $r'< l$.
	
	Since $\psi\in W^{1,\i}$, we have $u-\psi,\D(u-\psi)\in L^{r'}$ and $P(u-\psi)=0$. By assumption, the linear regularity statement gives $u-\psi\in W^{1,r'}$, and we may apply the $L^{r'}$-\emph{a priori} estimate~\eqref{L^p estimate} to $u-\psi\in\dom(\D_P)$:
	\begin{equation}
		\begin{split}\label{eq:similar_estimate_reg}
			\| u-\psi\|_{W^{1,r'}}^2 &\leq c_{r'} \left( \|u-\psi \|_{L^{r'}}^2 + \| \D  u-\D \psi \|_{L^{r'}}^2 \right) \\&\leq c_{r'} \left( (\|u\|_{L^{r'}} + \|\psi \|_{L^{r'}})^2 + (\| \D  u\|_{L^{r'}} + \|\D \psi \|_{L^{r'}})^2 \right).
		\end{split}
	\end{equation}
	Since $\psi\in W^{1,\i}(\Sigma \Omega)\subset W^{1,r'}(\Sigma \Omega)$, we get $\D \psi\in L^{r'}$ by Remark~\ref{remark: Lq_ineq}. This means that every term on the right-hand side of
	\[
	\| u\|_{W^{1,r'}} \leq \| \psi\|_{W^{1,r'}} + \sqrt{c_{r'}}\l( (\|u\|_{L^{r'}} + \|\psi \|_{L^{r'}})^2 + (\| \D  u\|_{L^{r'}} + \|\D \psi \|_{L^{r'}})^2 \r)^{\f12}
	\]
	has a finite upper bound, and $u\in W^{1,r'}(\Sigma \Omega)$.
	
	We now repeat the above inductively till we show that $u\in W^{1,\tilde r}$ for some $\tilde r\geq n$. For our choice of $p$, the generated exponents increase, since
	\[
	\frac{r'}r=\frac{n}{(n-r)(p-1)}>1
	\qquad\text{for }2\leq r<n.
	\]
	In particular, all exponents used in the regularity statement are greater than one.
	
	Assume $r'<n$ (otherwise we can skip the following steps), then $u \in L^{l'}$ with
	\begin{equation*}
		\f{1}{l'}= \frac{1}{r'} - \frac{1}{n} = \frac{p-1}{l} - \frac{1}{n}.
	\end{equation*}
	By the partial differential equation, we get $\D u \in L^{\frac{l'}{p-1}}(\Sigma \Omega) \allowbreak = L^{r''}(\Sigma \Omega)$ by setting $r'' \coloneq \frac{l'}{p-1}$.
	We get by the $L^{r''}$ regularity statement and \emph{a priori} estimate~\eqref{L^p estimate}, and the regularity of $\psi\in W^{1,\i}(\Sigma\Omega)\subset W^{1,r''}(\Sigma\Omega)$ as in~\eqref{eq:similar_estimate_reg}, that $u\in W^{1,r''}(\Sigma \Omega)$.
	If $r''<n$, we can characterize the next Sobolev exponent $l''$ by
	\begin{equation*}
		\f{1}{l''} = \frac{1}{r''} - \frac{1}{n} = \frac{p-1}{l'} - \frac{1}{n} =\frac{(p-1)^2}{l} - \frac{1+(p-1)}{n}.
	\end{equation*}
	
	Set $r^{(0)}=2$ and $l^{(0)}=l$.
	The $M\,$th step of this process assumes $u \in W^{1,r^{(M-1)}}(\Sigma \Omega)$ and $u\in L^{l^{(M-1)}}(\Sigma \Omega)$ with $r^{(M-1)}< n$. The partial differential equation implies $\D u\in L^{r^{(M)}}(\Sigma \Omega)$ with $r^{(M)}\coloneq \frac{l^{(M-1)}}{p-1}$, resulting in the regularity $u\in W^{1,r^{(M)}}(\Sigma \Omega)$ by reasoning as in~\eqref{eq:similar_estimate_reg} through the $L^{r^{(M)}}$ regularity statement and \emph{a priori} estimate~\eqref{L^p estimate} and the regularity of $\psi\in W^{1,\i}(\Sigma \Omega)\subset W^{1,r^{(M)}}(\Sigma \Omega)$.
	If $r^{(M)}\geq n$, we stop. Otherwise, the Sobolev embedding then implies that $u\in L^{l^{(M)}}(\Sigma \Omega)$ with
	\begin{equation}
		\frac{1}{l^{(M)}} = \frac{1}{r^{(M)}} - \frac{1}{n}= \frac{p-1}{l^{(M-1)}} - \frac{1}{n}.
		\label{eq:r^M}
	\end{equation}
	As long as the process continues, this inductive equation for $\frac{1}{l^{(M)}}$ can be explicitly calculated by
	\begin{equation}
		\begin{split}
			\frac{1}{l^{(M)}} &= \frac{(p-1)^M}{l} - \f{1}{n} \sum_{m=0}^{M-1} (p-1)^m \\
			&= \frac{(p-1)^M}{l} - \f{1}{n} \frac{1 - (p-1)^M}{1-(p-1)} \\
			&= (p-1)^M\bigg( \frac{1}{l}-\frac{1}{n(p-2)}\bigg) + \f{1}{n(p-2)}.
			\label{eq:recurrence}
		\end{split}
	\end{equation}
	Since $\frac{1}{l}-\frac{1}{n(p-2)}= \frac{n-2}{2n} - \frac{n-1}{2n} =-\frac{1}{2n}<0$, the last line \eqref{eq:recurrence} becomes negative for some sufficiently large $M$. If every $r^{(M)}$ remained below $n$, then \eqref{eq:r^M} would give $1/l^{(M)}>0$ for every $M$, contradicting \eqref{eq:recurrence}. Thus some $r^{(M)}\geq n$ was achieved, and the process stops at that step. Since $u\in W^{1,r^{(M)}}(\Sigma\Omega)$, we conclude that $u\in W^{1,n}(\Sigma\Omega)$.
	
	Having obtained $u\in W^{1,n}(\Sigma\Omega)$, the
	critical Sobolev embedding gives
	\[
	u\in L^b(\Sigma\Omega)
	\qquad\text{for every finite }b\geq1.
	\]
	Fix $1<a<\infty$ and choose $b=a(p-1)$.
	Then
	\[
	F(u)\in L^a,
	\qquad
	\|F(u)\|_{L^a}
	=
	\|u\|_{L^{a(p-1)}}^{p-1}.
	\]
	If the prescribed boundary datum admits an extension
	$\psi\in\bigcap_{1<a<\infty}W^{1,a}$, the section
	$v=u-\psi$ satisfies
	\[
	Pv=0,
	\qquad
	v,\D v\in L^a,
	\qquad
	\D v=\lambda F(u)-\D\psi.
	\]
	The assumed global Sobolev regularity therefore yields
	$v\in W^{1,a}$ and hence $u\in W^{1,a}$.
	By Morrey's inequality (Thm.~\ref{sobolev embeddings} (ii)), this implies $u \in C^{0,\alpha}(\Sigma \Omega)$ for all $\alpha\in (0,1)$.

	The nonlinearity $F(v)=|v|^{p-2}v$ is Lipschitz on bounded sets: its real derivative satisfies $DF(0)=0$ and
	$\|DF(v)\|\leq(p-1)|v|^{p-2}$ for $p>2$, so
	\[
	|F(v)-F(w)|\leq(p-1)R^{p-2}|v-w|
	\qquad\text{if }|v|,|w|\leq R.
	\]
	In smooth local orthonormal frames, this shows that
	$u\in C^{0,\alpha}$ implies
	$F(u)\in C^{0,\alpha}$, also at zeros.
	Furthermore, the H\"older regularity of $u$ implies the existence of a classical trace $\restr{u}{\partial\Omega}\in C^{0,\alpha}(\restr{\Sigma \Omega}{\partial \Omega})$.
	
	By interior Schauder regularity for Dirac operators \cite[Cor.~3.1.14]{Ammann03hab}, the $C^{0,\alpha}$ regularity of $|u|^{p-2}u$ in the equation $\D u = \la |u|^{p-2}u$ implies $u \in C^{1,\alpha}_{\mathrm{loc}}(\Sigma\operatorname{int}(\Omega))$ for all $\alpha \in (0,1)$. Under the additional boundary Schauder regularity assumption, and for a $C^{1,\alpha}$ extension $\psi$ of the prescribed boundary datum, we apply boundary regularity to
	\[
	P(u-\psi)=0,\qquad
	\D(u-\psi)=\la F(u)-\D\psi\in C^{0,\alpha},
	\]
	and obtain $u\in C^{1,\alpha}(\Sigma\Omega)$.
	
	Since $v\mapsto |v|^{p-2}v$ is smooth away from the origin, iterating higher interior Schauder estimates \cite[Thm.~3.1.16]{Ammann03hab} away from the zero set yields smoothness of $u$ in the interior. With smooth boundary data and the corresponding local boundary Schauder regularity of all orders, the same holds locally up to the boundary away from $u^{-1}(\{0\})$. The proof does not assert optimal regularity at the zero set.
	
	Lastly, note that for $n=2$, we can take $r=2=n$, implying by the critical Sobolev embedding that $u\in L^l$ for every finite $l$. For each $1<q<\infty$, choosing $l=q(p-1)$ gives $\D u\in L^q$. Applying the regularity statement to $u-\psi$ gives $u\in W^{1,q}$ for every finite $q>1$, and the conclusions follow as above. Only after Morrey's inequality gives boundedness of $u$ can we infer $\D u=\la|u|^{p-2}u\in L^\infty$.
	
	For the purely interior conclusions, the same exponent iteration uses the interior regularity statement of Remark~\ref{interior a priori} on nested relatively compact sets, without any boundary assumptions. Localization introduces only the zeroth-order term $[\D,\chi]u$, which has the required integrability at each step. Interior Morrey and Schauder regularity then give the stated conclusions.
\end{proof}

\begin{remark}[Other exponents and operators]
	For $n>2$, the same nonlinear bootstrap applies when
	\[
	2<p<\frac{2n-2}{n-2}.
	\]
	Indeed, for $2\leq r<n$,
	\[
	\frac{r'}{r}
	=
	\frac{n}{(n-r)(p-1)}
	\geq
	\frac{n}{(n-2)(p-1)}
	>1,
	\]
	so the iteration reaches $W^{1,n}$ after finitely
	many steps.
	For $n=2$, the argument applies to every finite $p>2$.
	
	The case $p=2$ is linear.
	The equation $\D u=\lambda u$ can be treated directly
	by elliptic regularity and yields the corresponding
	regularity conclusions under the same boundary
	hypotheses.
	
	At the endpoint
	$p=(2n-2)/(n-2)$ for $n>2$, the fixed-point existence
	argument remains valid, but the elementary bootstrap
	above gives no initial gain in integrability.
	No endpoint higher-regularity conclusion is asserted
	here without an additional argument.
	
	The reasoning also applies to other smooth first-order
	elliptic operators on vector bundles when the required
	Sobolev and Schauder regularity statements hold.
	On nonsmooth domains, regularity up to the boundary
	must be formulated separately; the fixed-point
	existence theorem alone does not imply such regularity.
\end{remark}

\section{Brezis--Nirenberg-type boundary value problem} \label{section:brezis--Nirenberg}
The same method applies to the spinorial Brezis--Nirenberg-type
boundary value problem
\[
\D u=a u+\lambda |u|^{p-2}u,
\qquad Pu=P\psi_0,
\qquad p=\frac{2n}{n-1}.
\]
Indeed, for
$a\in\mathbb R\setminus\operatorname{spec}(\D_P)$,
one replaces $\D_P$ by $\D_P-a$ and the harmonic
representative of the boundary datum by its
$(\D-a)$-harmonic representative.
The fixed-point argument then gives existence,
uniqueness in a prescribed ball, continuous dependence
on the boundary data, and geometric convergence under
the corresponding smallness condition; see
Corollary~\ref{cor:linear-perturbation}.
The regularity conclusions also remain valid.
As before, this is a small-data result for an
inhomogeneous boundary value problem, not an existence
result for nonzero solutions with homogeneous boundary
data.
\begin{corollary}[Application to the spinorial Brezis--Nirenberg problem]
	\label{cor:linear-perturbation}
	Let $\D_P$ be a self-adjoint and regular realization
	of $\D$ on a smooth compact Riemannian spin manifold
	of dimension $n\geq2$.
	Let
	\[
	a\in\mathbb R\setminus\operatorname{spec}(\D_P),
	\qquad
	p=\frac{2n}{n-1},
	\qquad
	q=2(p-1).
	\]
	Define
	\[
	A_q(a)=\|(\D_P-a)^{-1}\|_{\mathcal L(L^2,L^q)},
	\qquad
	B_1(a)=\|(\D_P-a)^{-1}\|_{\mathcal L(L^2,H^1)}.
	\]
	For a prescribed extension $\psi_0\in H^1(\Sigma\Omega)$,
	set
	\[
	\psi_a
	=
	\psi_0-(\D_P-a)^{-1}\bigl((\D-a)\psi_0\bigr).
	\]
	Then $\psi_a$ is the unique section in $H^1$ satisfying
	\[
	(\D-a)\psi_a=0,
	\qquad
	P\psi_a=P\psi_0.
	\]
	
	Let $\alpha>0$ and $M>0$ satisfy
	\[
	\|\psi_a\|_{L^q}\leq M,
	\qquad
	|\lambda|A_q(a)M^{p-2}<\delta_p(\alpha),
	\]
	where
	\[
	\delta_p(\alpha)
	=
	\min\left\{
	\frac{\alpha}{(1+\alpha)^{p-1}},
	\frac{1}{(p-1)(1+\alpha)^{p-2}}
	\right\}.
	\]
	Put
	\[
	R=\alpha M,
	\qquad
	L_M=(p-1)(1+\alpha)^{p-2}M^{p-2},
	\qquad
	\Theta_a=|\lambda|A_q(a)L_M<1.
	\]
	Then the following conclusions hold.
	
	\begin{enumerate}
		\item \emph{Existence and local uniqueness.}
		There exists a solution $u\in H^1(\Sigma\Omega)$ of
		\[
		\D u=a u+\lambda|u|^{p-2}u,
		\qquad
		Pu=P\psi_0,
		\]
		satisfying
		\[
		\|u-\psi_a\|_{L^q}\leq R.
		\]
		It is unique among all $H^1$-solutions satisfying
		this bound.
		
		\item \emph{Geometric convergence and estimates.}
		For every $u_0\in L^q(\Sigma\Omega)$ satisfying
		$\|u_0-\psi_a\|_{L^q}\leq R$, the iteration
		\[
		u_{k+1}
		=
		\psi_a+\lambda (\D_P-a)^{-1}\bigl(|u_k|^{p-2}u_k\bigr)
		\]
		converges to $u$ geometrically in $L^q$ and, after
		the first step, in $H^1$.
		More precisely, for $k\geq0$,
		\[
		\|u_k-u\|_{L^q}
		\leq
		\frac{\Theta_a^k}{1-\Theta_a}
		\|u_1-u_0\|_{L^q},
		\]
		and
		\[
		\|u_{k+1}-u\|_{H^1}
		\leq
		|\lambda|B_1(a)L_M\|u_k-u\|_{L^q}.
		\]
		Every iterate with $k\geq1$ belongs to $H^1$ and
		satisfies $Pu_k=P\psi_0$.
		Furthermore,
		\[
		\|u-\psi_a\|_{H^1}
		\leq
		|\lambda|B_1(a)(1+\alpha)^{p-1}M^{p-1},
		\]
		and the same bound holds for
		$\|u_k-\psi_a\|_{H^1}$ when $k\geq1$.
		
		\item \emph{Continuous dependence on the boundary data.}
		Fix $a,\lambda,\alpha,M$ as above.
		Let $\psi_{a,1},\psi_{a,2}$ be two
		$(\D-a)$-harmonic representatives with
		\[
		\|\psi_{a,j}\|_{L^q}\leq M,
		\]
		and let $u_j$ be their corresponding solutions in
		the balls
		$\|u_j-\psi_{a,j}\|_{L^q}\leq R$.
		Writing $w_j=u_j-\psi_{a,j}$, one has
		\[
		\|w_1-w_2\|_{L^q}
		\leq
		\frac{\Theta_a}{1-\Theta_a}
		\|\psi_{a,1}-\psi_{a,2}\|_{L^q},
		\]
		\[
		\|u_1-u_2\|_{L^q}
		\leq
		\frac{1}{1-\Theta_a}
		\|\psi_{a,1}-\psi_{a,2}\|_{L^q},
		\]
		and
		\[
		\|w_1-w_2\|_{H^1}
		\leq
		\frac{|\lambda|B_1(a)L_M}{1-\Theta_a}
		\|\psi_{a,1}-\psi_{a,2}\|_{L^q}.
		\]
		Consequently,
		\[
		\begin{aligned}
			\|u_1-u_2\|_{H^1}
			&\leq
			\|\psi_{a,1}-\psi_{a,2}\|_{H^1}\\
			&\quad+
			\frac{|\lambda|B_1(a)L_M}{1-\Theta_a}
			\|\psi_{a,1}-\psi_{a,2}\|_{L^q}.
		\end{aligned}
		\]
		
		\item \emph{Regularity.}
		The solution satisfies
		\[
		u\in C_{\mathrm{loc}}^{1,\gamma}
		\bigl(\Sigma(\operatorname{int}\Omega)\bigr)
		\qquad\text{for every }0<\gamma<1,
		\]
		and is smooth in the interior away from its zero set.
		
		If $(\D,P)$ satisfies the global Sobolev regularity
		assumptions of Theorem~\ref{second main theorem}
		and the prescribed boundary datum admits an extension
		\[
		\psi\in\bigcap_{1<r<\infty}W^{1,r}(\Sigma\Omega),
		\]
		then
		\[
		u\in W^{1,r}(\Sigma\Omega)
		\qquad\text{for every }1<r<\infty.
		\]
		For a fixed $0<\gamma<1$, the additional boundary
		Schauder regularity assumption and a
		$C^{1,\gamma}$ extension of the datum give
		$u\in C^{1,\gamma}(\Sigma\Omega)$.
		With a smooth extension and local boundary Schauder
		regularity of every order, $u$ is smooth away from
		its zero set, also locally up to the boundary.
	\end{enumerate}
	
	In particular, if $\psi_a\neq0$, taking
	$M=\|\psi_a\|_{L^q}$ and $\alpha=1/(p-2)$ gives
	the sufficient condition
	\[
	|\lambda|A_q(a)\|\psi_a\|_{L^q}^{p-2}
	<
	\frac{(p-2)^{p-2}}{(p-1)^{p-1}},
	\]
	with uniqueness in the ball
	\[
	\|u-\psi_a\|_{L^q}
	\leq
	\frac{\|\psi_a\|_{L^q}}{p-2}.
	\]
\end{corollary}

\begin{proof}
	Since $a\in\mathbb R$, the operator $-aI$ is a bounded
	self-adjoint operator on $L^2$. Thus $\D_P-a$ is
	self-adjoint on the unchanged domain
	\[
	\dom(\D_P-a)
	=
	\dom(\D_P)
	=
	\{v\in H^1(\Sigma\Omega):Pv=0\}.
	\]
	Moreover, the regularity estimate for $\D_P$ gives
	\[
	\begin{aligned}
		\|v\|_{H^1}^2
		&\leq
		c_1\bigl(\|v\|_{L^2}^2+\|\D_Pv\|_{L^2}^2\bigr)\\
		&\leq
		c_1\bigl((1+2a^2)\|v\|_{L^2}^2
		+2\|(\D_P-a)v\|_{L^2}^2\bigr).
	\end{aligned}
	\]
	Hence $\D_P-a$ is also regular.
	
	Set $
	\mu_a \coloneq \operatorname{dist}
	\bigl(a,\operatorname{spec}(\D_P)\bigr)>0.$
	The spectral theorem yields
	\[
	\|(\D_P-a)^{-1}f\|_{L^2}\leq\mu_a^{-1}\|f\|_{L^2}.
	\]
	Applying the preceding regularity estimate to $(\D_P-a)^{-1}f$
	therefore gives
	\[
	\|(\D_P-a)^{-1}f\|_{H^1}^2
	\leq
	c_1\bigl((1+2a^2)\mu_a^{-2}+2\bigr)
	\|f\|_{L^2}^2.
	\]
	Consequently, $B_1(a)<\infty$, and the Sobolev
	embedding $H^1\hookrightarrow L^q$ implies
	$A_q(a)<\infty$.
	
	\medskip
	\noindent\emph{The representative of the boundary datum.}
	Since $\psi_0\in H^1$, we have $(\D-a)\psi_0\in L^2$.
	Thus
	\[
	\psi_a=\psi_0-(\D_P-a)^{-1}((\D-a)\psi_0)
	\]
	belongs to $H^1$. The correction
	$(\D_P-a)^{-1}((\D-a)\psi_0)$ belongs to $\dom(\D_P)$ and
	therefore has homogeneous boundary data. It follows that
	\[
	P\psi_a=P\psi_0,
	\qquad
	(\D-a)\psi_a
	=
	(\D-a)\psi_0-(\D_P-a)(\D_P-a)^{-1}((\D-a)\psi_0)=0.
	\]
	If two sections have these properties, their difference
	belongs to $\dom(\D_P-a)$ and lies in its kernel.
	Invertibility therefore proves uniqueness.
	In particular, $\psi_a$ depends only on $P\psi_0$.
	
	\medskip
	\noindent\emph{The contraction argument.}
	Writing $u=\psi_a+w$, the boundary value problem is
	equivalent to
	\[
	(\D_P-a)w=\lambda F(w+\psi_a),
	\qquad
	F(v)=|v|^{p-2}v.
	\]
	Define
	\[
	T_a(w)=\lambda (\D_P-a)^{-1}F(w+\psi_a)
	\]
	on the closed ball
	\[
	\mathcal B_R
	=
	\{w\in L^q(\Sigma\Omega):\|w\|_{L^q}\leq R\},
	\qquad R=\alpha M.
	\]
	For $w\in\mathcal B_R$,
	\[
	\|w+\psi_a\|_{L^q}\leq R+M=(1+\alpha)M.
	\]
	Using the growth estimate for $F$, we obtain
	\[
	\|T_a(w)\|_{L^q}
	\leq
	|\lambda|A_q(a)(R+M)^{p-1}.
	\]
	The first term in the definition of $\delta_p(\alpha)$
	and the smallness assumption give
	\[
	\begin{aligned}
		|\lambda|A_q(a)(R+M)^{p-1}
		&=
		|\lambda|A_q(a)M^{p-2}
		(1+\alpha)^{p-1}M\\
		&<
		\alpha M=R.
	\end{aligned}
	\]
	Thus $T_a$ maps $\mathcal B_R$ into itself.
	
	For $w,z\in\mathcal B_R$, the Lipschitz estimate
	of Lemma~\ref{lem:F} gives
	\[
	\begin{aligned}
		\|T_a(w)-T_a(z)\|_{L^q}
		&\leq
		|\lambda|A_q(a)
		\|F(w+\psi_a)-F(z+\psi_a)\|_{L^2}\\
		&\leq
		|\lambda|A_q(a)(p-1)(R+M)^{p-2}
		\|w-z\|_{L^q}\\
		&=\Theta_a\|w-z\|_{L^q}.
	\end{aligned}
	\]
	The second term in $\delta_p(\alpha)$ ensures
	$\Theta_a<1$. Since $\mathcal B_R$ is complete,
	the Banach fixed-point theorem gives a unique fixed
	point $w\in\mathcal B_R$.
	
	The fixed-point identity implies
	\[
	w=\lambda (\D_P-a)^{-1}F(w+\psi_a)\in\dom(\D_P-a).
	\]
	Hence $u=w+\psi_a$ belongs to $H^1$ and satisfies
	both the equation and the prescribed boundary condition.
	Conversely, every $H^1$-solution in the stated ball
	gives a fixed point of $T_a$, proving the asserted
	local uniqueness.
	
	\medskip
	\noindent\emph{Convergence and estimates.}
	For $w_k=u_k-\psi_a$, the iteration is
	$w_{k+1}=T_a(w_k)$. The contraction theorem yields
	\[
	\|w_k-w\|_{L^q}
	\leq
	\frac{\Theta_a^k}{1-\Theta_a}
	\|w_1-w_0\|_{L^q}.
	\]
	Moreover,
	\[
	w_{k+1}-w
	=
	\lambda (\D_P-a)^{-1}\bigl(F(u_k)-F(u)\bigr).
	\]
	Since $\|u_k\|_{L^q},\|u\|_{L^q}\leq R+M$,
	\[
	\|w_{k+1}-w\|_{H^1}
	\leq
	|\lambda|B_1(a)L_M\|u_k-u\|_{L^q}.
	\]
	This proves geometric convergence in $H^1$ after
	the first step. Similarly,
	\[
	\|w_{k+1}\|_{H^1},\ \|w\|_{H^1}
	\leq
	|\lambda|B_1(a)(R+M)^{p-1},
	\]
	which gives the stated uniform estimates.
	Each $w_{k+1}$ belongs to $\dom(\D_P)$, so every
	$u_{k+1}$ satisfies the prescribed boundary condition.
	
	\medskip
	\noindent\emph{Continuous dependence.}
	For the two representatives and corresponding solutions,
	write $w_j=u_j-\psi_{a,j}$. Their fixed-point identities
	give
	\[
	w_1-w_2=\lambda (\D_P-a)^{-1}\bigl(F(u_1)-F(u_2)\bigr).
	\]
	Since $\|u_j\|_{L^q}\leq R+M$, it follows that
	\[
	\|w_1-w_2\|_{L^q}
	\leq
	\Theta_a\|u_1-u_2\|_{L^q}.
	\]
	Combining this with
	\[
	\|u_1-u_2\|_{L^q}
	\leq
	\|w_1-w_2\|_{L^q}
	+
	\|\psi_{a,1}-\psi_{a,2}\|_{L^q}
	\]
	and absorbing the term with coefficient $\Theta_a$
	yields
	\[
	\|u_1-u_2\|_{L^q}
	\leq
	\frac{1}{1-\Theta_a}
	\|\psi_{a,1}-\psi_{a,2}\|_{L^q}.
	\]
	Substitution gives the claimed estimate for
	$\|w_1-w_2\|_{L^q}$. Finally,
	\[
	\begin{aligned}
		\|w_1-w_2\|_{H^1}
		&\leq
		|\lambda|B_1(a)L_M\|u_1-u_2\|_{L^q}\\
		&\leq
		\frac{|\lambda|B_1(a)L_M}{1-\Theta_a}
		\|\psi_{a,1}-\psi_{a,2}\|_{L^q}.
	\end{aligned}
	\]
	The $H^1$ estimate for $u_1-u_2$ follows by the
	triangle inequality.
	
	\medskip
	\noindent\emph{Regularity.}
	The additional term $au$ does not alter the exponent
	iteration in the proof of
	Theorem~\ref{second main theorem}.
	Indeed, if $u\in L^\ell$ and
	$r=\ell/(p-1)>1$, then
	\[
	F(u)\in L^r,
	\qquad
	\|F(u)\|_{L^r}=\|u\|_{L^\ell}^{p-1}.
	\]
	Since $r<\ell$ and $\Omega$ has finite volume,
	also $au\in L^r$, with
	\[
	\|au\|_{L^r}
	\leq
	|a|\,\operatorname{vol}_g(\Omega)^{1/r-1/\ell}
	\|u\|_{L^\ell}.
	\]
	Thus $\D u=au+\lambda F(u)\in L^r$, and the
	same elliptic estimates and Sobolev embeddings
	give the same improvement of integrability.
	
	For the global argument, choose an extension
	$\psi\in\bigcap_{1<r<\infty}W^{1,r}(\Sigma\Omega)$ of the prescribed
	boundary datum and apply the homogeneous-boundary
	estimates to $u-\psi$:
	\[
	P(u-\psi)=0,
	\qquad
	\D(u-\psi)=au+\lambda F(u)-\D\psi.
	\]
	The term $\D\psi$ belongs to every finite $L^r$,
	so the bootstrap proceeds as before.
	The interior argument uses local elliptic estimates
	and does not require any boundary regularity assumption.
	
	Once $u\in C^{0,\gamma}$, the local Lipschitz
	continuity of $F$ implies
	$au+\lambda F(u)\in C^{0,\gamma}$.
	Schauder estimates therefore give $C^{1,\gamma}$
	regularity, in the interior and, under the stated
	additional assumptions, up to the boundary.
	Away from the zero set, the right-hand side is a
	smooth function of $u$, and iteration of higher
	Schauder estimates gives smoothness.
	
	Finally, choosing $M=\|\psi_a\|_{L^q}>0$ and
	$\alpha=1/(p-2)$ gives the optimized smallness
	condition and uniqueness radius stated above.
\end{proof}

\begin{remark}[Homogenous Problem] \label{rem:homogeneous-gap-Brezis}
	If $P\psi_0=0$ (or even if $\Omega$ is closed), then $\psi_a=0$, and the contraction
	constructs only the zero solution.
	For $\lambda\neq0$, every nonzero $H^1$-solution
	with homogeneous boundary data satisfies
	\[
	\|u\|_{L^q}
	\geq
	\bigl(|\lambda|A_q(a)\bigr)^{-1/(p-2)}.
	\]
	This is analogous to Remark~\ref{rem:closed} and Corollary~\ref{cor:Lq-gap}.
\end{remark}

\section{Appendix}\label{appendix}
We will show that the Euler--Lagrange equation associated
to the variational characterization of
$\lambda_{\mathrm{min}}(\Omega,\partial\Omega)$ is
equivalent, after normalization, to the spinorial Yamabe
equation~\eqref{def:pde_homogenous}.
Throughout, we assume that $\D_B$ is self-adjoint and
invertible, with
\[
\dom(\D_B)
=
\{\varphi\in H^1(\Sigma\Omega):B\varphi=0\}.
\]

\subsection{The Variational Characterization}
The reason why specifically the chiral boundary
condition $B^\pm$ appears in the definition of the
conformal invariant is mainly because chirality
operators behave nicely under conformal change of
the metric $g$ \cite{Raulot08}.
Using these transformation behaviours, Raulot obtained
the corresponding variational characterization.
In the invertible setting considered here, we write it as
\[
\lambda_{\mathrm{min}}(\Omega,\partial\Omega)
=
\inf_{g'\in[g]}
|\lambda_1(\D_B^{g'})|
\,\vol(\Omega,g')^{\frac1n}
=
\inf_{\varphi\in\mathcal U}
\frac{
	\left(\int_\Omega|\D\varphi|^q\,d\vol_g\right)^{2/q}}
{\left|
	\int_\Omega
	\operatorname{Re}\langle\D\varphi,\varphi\rangle
	\,d\vol_g
	\right|},
\]
where
\[
q=\frac{2n}{n+1},
\qquad
p=\frac{2n}{n-1},
\qquad
\frac1p+\frac1q=1,
\]
and
\[
\mathcal U
=
\{\varphi\in\dom(\D_B):J(\varphi)\neq0\},
\qquad
J(\varphi)
=
\int_\Omega
\operatorname{Re}\langle\D\varphi,\varphi\rangle
\,d\vol_g.
\]
The restriction $J(\varphi)\neq0$ is necessary:
invertibility of $\D_B$ does not prevent cancellation
between positive and negative spectral contributions.

Set
\[
S(\varphi)
=
\int_\Omega|\D\varphi|^q\,d\vol_g.
\]
Taking the quotient above as our functional gives
\[
\mathcal F:\mathcal U\longrightarrow\mathbb R_{>0},
\qquad
\mathcal F(\varphi)
=
\frac{S(\varphi)^{2/q}}{|J(\varphi)|}.
\]
We regard $\dom(\D_B)$ as a real Banach space with
its graph norm. Since $q>1$, this functional is
continuously Fr\'echet differentiable on $\mathcal U$.

For $\eta\in\dom(\D_B)$,
\[
S'(\varphi)[\eta]
=
q\int_\Omega
\operatorname{Re}
\left\langle
|\D\varphi|^{q-2}\D\varphi,\D\eta
\right\rangle\,d\vol_g.
\]
Self-adjointness gives
\[
J'(\varphi)[\eta]
=
2\int_\Omega
\operatorname{Re}\langle\varphi,\D\eta\rangle
\,d\vol_g.
\]
Consequently, $\mathcal F'(\varphi)=0$ is equivalent
to
\[
\int_\Omega
\operatorname{Re}
\left\langle
|\D\varphi|^{q-2}\D\varphi
-\frac{S(\varphi)}{J(\varphi)}\varphi,
\D\eta
\right\rangle\,d\vol_g
=
0
\]
for every $\eta\in\dom(\D_B)$.
Equivalently, using the shorthand
\[
A[\varphi]
=
2\frac{J(\varphi)}{S(\varphi)},
\]
this is
\[
\int_\Omega
\operatorname{Re}
\left\langle
A[\varphi]|\D\varphi|^{q-2}\D\varphi
-2\varphi,\D\eta
\right\rangle\,d\vol_g
=
0.
\]

Define
\[
\Psi\coloneq|\D\varphi|^{q-2}\D\varphi,
\]
where this expression is defined to be zero when
$\D\varphi=0$.
Since $(q-1)p=q$, we have
\[
\int_\Omega|\Psi|^p\,d\vol_g
=
\int_\Omega|\D\varphi|^q\,d\vol_g
=
S(\varphi).
\]
Thus $\Psi\in L^p(\Sigma\Omega)\subset L^2(\Sigma\Omega)$.

Now use invertibility of $\D_B$: its range is all of
$L^2(\Sigma\Omega)$. The critical-point identity
therefore implies
\[
\Psi=\lambda\varphi,
\qquad
\lambda\coloneq
\frac{S(\varphi)}{J(\varphi)}
=
\frac{2}{A[\varphi]}
\in\mathbb R\setminus\{0\}.
\]
In particular,
\[
\Psi\in\dom(\D_B),
\qquad
B\Psi=0.
\]
This establishes the boundary condition without
integrating by parts on $\Psi$ before its regularity
is known. Applying $\D_B$ now gives
\[
\D\Psi=\lambda\D\varphi.
\]

Since $\mathcal F$ is positively $0$-homogeneous,
\[
\mathcal F(t\varphi)=\mathcal F(\varphi)
\qquad(t>0),
\]
we may normalize a critical point by
\[
\|\D\varphi\|_{L^q}=1.
\]
The Euler--Lagrange equation can consequently be
written as
\begin{equation}
	\begin{gathered}
		\D\bigl(|\D\varphi|^{q-2}\D\varphi\bigr)
		=
		\lambda\D\varphi,
		\qquad
		\|\D\varphi\|_{L^q}=1,
		\qquad
		B\varphi=0,\\
		|\D\varphi|^{q-2}\D\varphi\in\dom(\D_B),
		\qquad
		\lambda=\frac1{J(\varphi)}.
	\end{gathered}
	\label{eq:EulerLagrange_q}
\end{equation}
The domain condition on the nonlinear spinor is part
of this boundary Euler--Lagrange problem, not merely
a consequence of its interior differential equation.

This Euler--Lagrange equation is equivalent to the
normalized spinorial Yamabe equation in the following
way. The maps
\[
z\longmapsto |z|^{q-2}z,
\qquad
z\longmapsto |z|^{p-2}z
\]
are inverses, since $(p-1)(q-1)=1$.
Thus, with
$\Psi=|\D\varphi|^{q-2}\D\varphi$, we have
\[
\D\varphi=|\Psi|^{p-2}\Psi.
\]
Equation~\eqref{eq:EulerLagrange_q} becomes
\begin{equation}
	\D\Psi
	=
	\lambda|\Psi|^{p-2}\Psi,
	\qquad
	\|\Psi\|_{L^p}=1,
	\qquad
	B\Psi=0.
	\label{eq:EulerLagrange_p}
\end{equation}
Moreover, the relation between the two spinors is
$\Psi=\lambda\varphi$.

Conversely, let $\Psi\in\dom(\D_B)$ solve
\eqref{eq:EulerLagrange_p}, with
$\lambda\in\mathbb R\setminus\{0\}$, and put
\[
\varphi\coloneq\frac{\Psi}{\lambda}.
\]
Then $\varphi\in\dom(\D_B)$ and
\[
\D\varphi=|\Psi|^{p-2}\Psi,
\qquad
|\D\varphi|^{q-2}\D\varphi=\Psi.
\]
Since $(p-1)q=p$, we obtain
\[
S(\varphi)
=
\int_\Omega|\D\varphi|^q\,d\vol_g
=
\int_\Omega|\Psi|^p\,d\vol_g
=
1.
\]
Furthermore,
\[
J(\varphi)
=
\frac1\lambda
\int_\Omega|\Psi|^p\,d\vol_g
=
\frac1\lambda.
\]
Hence $\varphi\in\mathcal U$ and
$S(\varphi)/J(\varphi)=\lambda$.
The identity
\[
|\D\varphi|^{q-2}\D\varphi
=
\lambda\varphi
\]
therefore gives $\mathcal F'(\varphi)=0$.
This proves the correspondence between normalized
critical points and normalized $H^1$ solutions
of the spinorial Yamabe equation.

At such a normalized critical point,
\[
\mathcal F(\varphi)=|\lambda|.
\]
In particular, if the infimum is attained by a
critical point in this class, then its associated
solution satisfies
$|\lambda|=\lambda_{\mathrm{min}}(\Omega,\partial\Omega)$.
The correspondence itself does not establish
attainment of the infimum.

For another boundary operator $P$ with an invertible
self-adjoint realization
\[
\dom(\D_P)
=
\{v\in H^1(\Sigma\Omega):Pv=0\},
\]
the same algebraic and variational argument applies
to the quotient
\[
\inf_{\substack{\varphi\in\dom(\D_P)\\J(\varphi)\neq0}}
\frac{
	\left(\int_\Omega|\D\varphi|^q\,d\vol_g\right)^{2/q}}
{|J(\varphi)|}.
\]
However, identifying this quantity with a conformal
spectral invariant requires additional conformal
covariance properties of the boundary condition.
	
\subsection{Application to small geodesic balls}
\label{subsec:geodesic-balls}

To apply our existence theorem on a geodesic ball
$\Omega_r=B_r^g(x_0)$, we need an upper bound for
\[
A_q(\Omega_r)
=
\|\D_{P_r}^{-1}\|_{\mathcal L(L^2(\Omega_r),L^q(\Omega_r))}.
\]
Rather than estimating the Gagliardo--Nirenberg constant,
the regularity constant, and the spectral gap separately,
we compare the inverse Dirac operator with its Euclidean
counterpart.

We first rescale the geodesic ball to the unit ball.
The rescaled metric approaches the Euclidean metric,
and the Bourguignon--Gauduchon identification allows
us to compare the corresponding Dirac operators as
maps from $H^1$ to $L^2$.
We then use a Neumann-series argument to estimate
the inverse and finally undo the rescaling.

Throughout this subsection, we assume that the spinor
bundle carries a chirality operator.
In odd dimension, we use the doubled spinor bundle
introduced in the section on chiral boundary conditions.

\begin{theorem}
	\label{thm:geodesic-ball-inverse}
	Let $(M,g)$ be a smooth Riemannian spin manifold of
	dimension $n\geq2$, and fix an interior point $x_0\in M$.
	Set
	\[
	\Omega_r=B_r^g(x_0),
	\qquad
	q=\frac{2(n+1)}{n-1}.
	\]
	For each sufficiently small $r>0$, let
	\[
	P_r^\pm u
	=
	\frac12\bigl(\id\mp i\gamma_g(\nu_r)\beta\bigr)
	u|_{\partial\Omega_r},
	\]
	where $\nu_r$ is the outward unit normal to
	$\partial\Omega_r$.
	Fix one of the two signs and write $P_r$ for the
	corresponding boundary operator.
	
	There exist $r_0>0$ and $C>0$ such that, for every
	$0<r\leq r_0$, the realization $\D_{P_r}$ is invertible
	and
	\[
	A_q(\Omega_r)\leq Cr^{1/(n+1)}.
	\]
\end{theorem}

\begin{proof}
	Choose an oriented isometry $T_{x_0}M\cong\mathbb R^n$.
	For sufficiently small $r>0$, the map
	\[
	\Phi_r:B_1(0)\longrightarrow\Omega_r,
	\qquad
	\Phi_r(y)=\exp_{x_0}(ry),
	\]
	is a diffeomorphism.
	Define the metric
	\[
	g_r=r^{-2}\Phi_r^*g
	\]
	on $B_1(0)$.
	In normal coordinates, its coefficients are
	$(g_r)_{ij}(y)=g_{ij}(ry)$.
	The normal-coordinate expansion therefore gives
	\begin{equation}
		\|g_r-g_{\mathrm{eucl}}\|_{C^1(\overline{B_1})}
		\leq Cr^2.
		\label{eq:ball-metric-comparison}
	\end{equation}
	
	Let $G_r$ be the positive symmetric endomorphism of
	$TB_1$ determined by
	\[
	g_r(X,Y)=g_{\mathrm{eucl}}(G_rX,Y),
	\]
	and put $
	B_r=G_r^{-1/2}.$
	Then
	\[
	g_r(B_rX,B_rY)=g_{\mathrm{eucl}}(X,Y).
	\]
	Thus $B_r$ maps positively oriented Euclidean
	orthonormal frames to positively oriented
	$g_r$-orthonormal frames.
	Lifting this map to the spin structures gives the
	Bourguignon--Gauduchon identification
	\[
	\xi_r:\Sigma_{\mathrm{eucl}}B_1
	\longrightarrow\Sigma_{g_r}B_1
	\]
	\cite{Bourguignon92,Hermann12}.
	It is fiberwise unitary and satisfies
	\begin{equation}
		\xi_r\bigl(\gamma_{\mathrm{eucl}}(X)v\bigr)
		=
		\gamma_{g_r}(B_rX)\xi_rv.
		\label{eq:ball-Clifford-identification}
	\end{equation}
	On doubled bundles, we use the componentwise
	identification.
	
	We use $\xi_r$ to compare the two Dirac operators
	on Euclidean spinors.
	Let $e_j=\partial_j$ be the Euclidean coordinate frame
	and write
	\[
	e_{r,j}=B_re_j,
	\qquad
	B_re_j=\sum_{k=1}^n b_{r,j}^{\,k}\partial_k.
	\]
	The frame $(e_{r,j})$ is orthonormal for $g_r$.
	In the corresponding spin frame, the spin connection
	has the form
	\[
	\xi_r^{-1}\nabla^{g_r}_{e_{r,j}}(\xi_rv)
	=
	\sum_{k=1}^n b_{r,j}^{\,k}\partial_kv
	+
	\Gamma_{r,j}v,
	\]
	where $\Gamma_{r,j}$ is a matrix-valued function.
	Consequently,
	\[
	\xi_r^{-1}\D^{g_r}\xi_r v
	=
	-i\sum_{j,k=1}^n
	\gamma_{\mathrm{eucl}}(e_j)b_{r,j}^{\,k}\partial_kv
	-i\sum_{j=1}^n
	\gamma_{\mathrm{eucl}}(e_j)\Gamma_{r,j}v.
	\]
	
	By \eqref{eq:ball-metric-comparison},
	\[
	\|B_r-I\|_{C^1(\overline{B_1})}\leq Cr^2.
	\]
	The Christoffel symbols of $g_r$ are also $O(r^2)$.
	Since the spin connection coefficients in the frame
	$(e_{r,j})$ depend on these symbols and the first
	derivatives of $B_r$, we obtain
	\[
	\|\Gamma_{r,j}\|_{L^\infty(B_1)}\leq Cr^2.
	\]
	Subtracting
	\[
	\D^{\mathrm{eucl}}v
	=
	-i\sum_{j=1}^n
	\gamma_{\mathrm{eucl}}(e_j)\partial_jv
	\]
	therefore gives
	\[
	\left|
	\left(\xi_r^{-1}\D^{g_r}\xi_r
	-\D^{\mathrm{eucl}}\right)v
	\right|
	\leq Cr^2\bigl(|\nabla^{\mathrm{eucl}}v|+|v|\bigr).
	\]
	After integration, there is
	a constant $C_E>0$, independent of sufficiently small
	$r$, such that
	\begin{equation}
		\left\|
		\left(\xi_r^{-1}\D^{g_r}\xi_r
		-\D^{\mathrm{eucl}}\right)v
		\right\|_{L^2}
		\leq C_Er^2\|v\|_{H^1},
		\qquad v\in H^1(\Sigma_{\mathrm{eucl}}B_1).
		\label{eq:ball-Dirac-comparison}
	\end{equation}
	In particular, the operator difference is small in
	$\mathcal L(H^1,L^2)$.
	
	To compare the inverses, we must also verify that
	$\xi_r$ maps the homogeneous Euclidean boundary domain
	onto the homogeneous boundary domain for $g_r$.
	
	Let $\beta_{\mathrm{eucl}}$ and $\beta_{g_r}$ denote
	chirality operators for the two metrics:
	In even dimension, they are given in the above
	positively oriented orthonormal frames by
	\[
	\beta_{\mathrm{eucl}}
	=
	i^{n/2}\gamma_{\mathrm{eucl}}(e_1)\cdots
	\gamma_{\mathrm{eucl}}(e_n),
	\qquad
	\beta_{g_r}
	=
	i^{n/2}\gamma_{g_r}(e_{r,1})\cdots
	\gamma_{g_r}(e_{r,n}).
	\]
	Equation~\eqref{eq:ball-Clifford-identification} yields
	\[
	\beta_{g_r}\xi_r=\xi_r\beta_{\mathrm{eucl}}.
	\]
	On the doubled bundles in odd dimension, use
	$\beta=\operatorname{diag}(\id,-\id)$ for both metrics.
	Since $\xi_r$ acts componentwise, it also intertwines
	these chirality operators.
	
	For $y\in\overline{B_1}$, let
	\[
	Y_y=\sum_{j=1}^n y^j
	\left.\frac{\partial}{\partial y^j}\right|_y
	\in T_yB_1
	\]
	be the radial vector.
	The Gauss lemma gives
	\[
	(g_r)_y(Y_y,Z)
	=
	(g_{\mathrm{eucl}})_y(Y_y,Z)
	\qquad
	\text{for every }Z\in T_yB_1.
	\]
	By the definition of $G_r$, this implies $
	G_r(y)Y_y=Y_y.$
	Since $B_r=G_r^{-1/2}$, it follows that $
	B_r(y)Y_y=Y_y. $
	
	For $y\in\partial B_1$, the vector $Y_y$ is the outward
	unit normal for both metrics.
	Consequently,
	\[
	B_r\nu_{\mathrm{eucl}}=\nu_{g_r},
	\]
	and the Clifford-intertwining identity gives
	\[
	\gamma_{g_r}(\nu_{g_r})\xi_r
	=
	\xi_r\gamma_{\mathrm{eucl}}(\nu_{\mathrm{eucl}}).
	\]
	Together with the corresponding identity for chirality,
	this yields
	\[
	P_{g_r}^\pm(\xi_rv)
	=
	\xi_r(P_{\mathrm{eucl}}^\pm v)
	\quad\text{on }\partial B_1,
	\]
	where the equality is understood in the trace sense
	for $v\in H^1$.
	
	The smooth bundle map $\xi_r$ and its inverse induce
	bounded isomorphisms between the corresponding $H^1$
	spaces.
	Indeed, in local trivializations, multiplication by
	their smooth matrix representatives preserves $H^1$,
	and the coordinate and intrinsic Sobolev norms are
	equivalent.
	Thus the transported realization
	$\xi_r^{-1}\D^{g_r}_{P_{g_r}^\pm}\xi_r$
	and the Euclidean chiral realization have the common
	domain
	\[
	\left\{
	v\in H^1(\Sigma_{\mathrm{eucl}}B_1):
	P_{\mathrm{eucl}}^\pm v=0
	\right\}.
	\]
	
	Let $\D_P^{\mathrm{eucl}}$ denote the Euclidean
	realization with the chosen chiral boundary condition.
	It is self-adjoint and regular.
	We first verify that it is invertible.
	
	We verify that $\D_P^{\mathrm{eucl}}$ has trivial kernel.
	In this computation, the Hermitian product is taken
	to be linear in its first argument.
	
	For $u,z\in H^1(\Sigma_{\mathrm{eucl}}B_1)$,
	Green's identity for
	$\D^{\mathrm{eucl}}=-i\sum_j\gamma_{\mathrm{eucl}}(e_j)
	\partial_j$ is
	\[
	\langle\D^{\mathrm{eucl}}u,z\rangle_{L^2}
	-
	\langle u,\D^{\mathrm{eucl}}z\rangle_{L^2}
	=
	-i\int_{\partial B_1}
	\langle\gamma_{\mathrm{eucl}}(\nu)u,z\rangle
	\,d\sigma.
	\]
	This follows by integration by parts for smooth
	sections and extends to $H^1$ by density and
	continuity of the trace.
	
	Write the chosen boundary projector as
	\[
	P^\varepsilon
	=
	\frac12
	\bigl(
	\id-\varepsilon i\gamma_{\mathrm{eucl}}(\nu)
	\beta_{\mathrm{eucl}}
	\bigr),
	\qquad
	\varepsilon\in\{1,-1\}.
	\]
	Let $v\in\ker(\D_{P^\varepsilon}^{\mathrm{eucl}})$.
	Since $\beta_{\mathrm{eucl}}$ is parallel and
	anticommutes with Clifford multiplication,
	\[
	\D^{\mathrm{eucl}}(\beta_{\mathrm{eucl}}v)
	=
	-\beta_{\mathrm{eucl}}\D^{\mathrm{eucl}}v=0.
	\]
	Applying Green's identity to $u=v$ and
	$z=\beta_{\mathrm{eucl}}v$ therefore gives
	\[
	0=
	-i\int_{\partial B_1}
	\langle
	\gamma_{\mathrm{eucl}}(\nu)v,
	\beta_{\mathrm{eucl}}v
	\rangle\,d\sigma.
	\]
	
	The boundary condition $P^\varepsilon v=0$ implies
	\[
	\gamma_{\mathrm{eucl}}(\nu)v
	=
	\varepsilon i\beta_{\mathrm{eucl}}v.
	\]
	Since $\beta_{\mathrm{eucl}}$ is unitary, we obtain
	\[
	0
	=
	\varepsilon
	\int_{\partial B_1}
	|\operatorname{Tr}v|^2\,d\sigma.
	\]
	Hence $\operatorname{Tr}v=0$.
	This conclusion uses both the equation
	$\D^{\mathrm{eucl}}v=0$ and the chiral boundary
	condition, not the boundary condition alone.
	
	Thus $v\in H^1_0(\Sigma_{\mathrm{eucl}}B_1)$.
	For every $z\in H^1_0(\Sigma_{\mathrm{eucl}}B_1)$,
	\[
	\|\D^{\mathrm{eucl}}z\|_{L^2}^2
	=
	\|\nabla^{\mathrm{eucl}}z\|_{L^2}^2.
	\]
	Indeed, for compactly supported smooth sections this
	follows from
	$(\D^{\mathrm{eucl}})^2=-\Delta$ and integration
	by parts, and the identity extends to $H^1_0$
	by density.
	Applying it to $v$ yields
	$\|\nabla^{\mathrm{eucl}}v\|_{L^2}=0$.
	Poincar\'e's inequality then gives $v=0$.
	
	The realization is self-adjoint with compact resolvent.
	Its trivial kernel therefore implies that it is
	invertible.
	
	Denote its inverse constants by
	\[
	A_q^{\mathrm{eucl}}
	=
	\|(\D_P^{\mathrm{eucl}})^{-1}\|
	_{\mathcal L(L^2,L^q)},
	\qquad
	B_1^{\mathrm{eucl}}
	=
	\|(\D_P^{\mathrm{eucl}})^{-1}\|
	_{\mathcal L(L^2,H^1)}.
	\]
	Both constants are finite by regularity and Sobolev
	embedding.

	Set
	\[
	E_r = \xi_r^{-1}\D^{g_r}\xi_r - \D^{\mathrm{eucl}}.
	\]
	By \eqref{eq:ball-Dirac-comparison}, $E_r$ extends to a bounded operator from $H^1(\Sigma_{\mathrm{eucl}}B_1)$ to $L^2(\Sigma_{\mathrm{eucl}}B_1)$ with $\|E_r\|_{\mathcal{L}(H^1, L^2)} \leq C_E r^2$. Since $(\D_P^{\mathrm{eucl}})^{-1}$ is a bounded operator from $L^2$ to $H^1$ with norm $B_1^{\mathrm{eucl}}$, the composition $E_r (\D_P^{\mathrm{eucl}})^{-1}$ is a bounded operator on $L^2(\Sigma_{\mathrm{eucl}}B_1)$ satisfying
	\[
	\|E_r (\D_P^{\mathrm{eucl}})^{-1}\|_{\mathcal{L}(L^2, L^2)} \leq C_E r^2 B_1^{\mathrm{eucl}}.
	\]
	Choose $r_0 > 0$ sufficiently small such that $C_E r_0^2 B_1^{\mathrm{eucl}} \leq 1/2$. Then, for all $0 < r \leq r_0$, the operator $I + E_r (\D_P^{\mathrm{eucl}})^{-1}$ is invertible on $L^2$ via a convergent Neumann series, and its inverse satisfies
	\begin{equation}
		\left\| \left( I + E_r (\D_P^{\mathrm{eucl}})^{-1} \right)^{-1} \right\|_{\mathcal{L}(L^2, L^2)} \leq \frac{1}{1 - C_E r^2 B_1^{\mathrm{eucl}}} \leq 2.
		\label{eq:ball-Neumann-estimate}
	\end{equation}
	
	Since the boundary domains agree, we have the operator factorization
	\[
	\xi_r^{-1}\D^{g_r}_{P_{g_r}^\pm}\xi_r = (I + E_r (\D_P^{\mathrm{eucl}})^{-1}) \D_P^{\mathrm{eucl}}.
	\]
	It follows that the transported realization is invertible, with inverse
	\[
	(\xi_r^{-1}\D^{g_r}_{P_{g_r}^\pm}\xi_r)^{-1} = (\D_P^{\mathrm{eucl}})^{-1} (I + E_r (\D_P^{\mathrm{eucl}})^{-1})^{-1}.
	\]
	Using the boundedness of $(\D_P^{\mathrm{eucl}})^{-1}$ and \eqref{eq:ball-Neumann-estimate}, we obtain the following bounds on the operator norms:
	\begin{align}
		A_q(B_1, g_r) &\coloneq \|(\D^{g_r}_{P_{g_r}^\pm})^{-1}\|_{\mathcal{L}(L^2, L^q)} \leq 2 A_q^{\mathrm{eucl}}, \label{eq:ball-reference-inverse-estimate} \\
		B_1(B_1, g_r) &\coloneq \|(\D^{g_r}_{P_{g_r}^\pm})^{-1}\|_{\mathcal{L}(L^2, H^1)} \leq 2 B_1^{\mathrm{eucl}}, \label{eq:ball-reference-H1-estimate}
	\end{align}
	where the norms are compared using the fact that $\xi_r$ preserves pointwise spinor lengths and the volume forms $d\operatorname{vol}_{g_r}$ and $d\operatorname{vol}_{\mathrm{eucl}}$ are equivalent with comparison factors tending to one as $r \to 0$.
	
	Finally, $\Phi_r$ identifies $(\Omega_r, g)$ with $(B_1, r^2 g_r)$. Under the canonical identification of spinors for constant metric scaling, the inverse operators satisfy $(\D^{r^2 g_r}_P)^{-1} = r (\D^{g_r}_P)^{-1}$. Consequently, for any $f \in L^2(\Sigma \Omega_r)$, the scaling of the norms yields
	\[
	\frac{\|(\D^{r^2 g_r}_P)^{-1} f\|_{L^q(r^2 g_r)}}{\|f\|_{L^2(r^2 g_r)}} = r^{1 + n/q - n/2} \frac{\|(\D^{g_r}_P)^{-1} \tilde{f}\|_{L^q(g_r)}}{\|\tilde{f}\|_{L^2(g_r)}},
	\]
	where $\tilde{f}$ is the pushed-forward section on $B_1$. Taking the supremum over $f \neq 0$ gives
	\begin{equation*}
		A_q(\Omega_r) = r^{1 + n/q - n/2} A_q(B_1, g_r).
	\end{equation*}
	Substituting $q = 2(n+1)/(n-1)$, the exponent becomes $1 + n/q - n/2 = 1/(n+1)$. Combined with \eqref{eq:ball-reference-inverse-estimate}, we conclude
	\[
	A_q(\Omega_r) \leq C r^{1/(n+1)}.
	\]
	\end{proof}
	
	\begin{remark}[Comparison with the Euclidean inverse]
	The perturbation argument also yields convergence of the inverse operators. The identity
	\[
	(\xi_r^{-1}\D^{g_r}_{P_{g_r}^\pm}\xi_r)^{-1} - (\D_P^{\mathrm{eucl}})^{-1} = - (\D_P^{\mathrm{eucl}})^{-1} E_r (\xi_r^{-1}\D^{g_r}_{P_{g_r}^\pm}\xi_r)^{-1}
	\]
	combined with \eqref{eq:ball-Dirac-comparison} and \eqref{eq:ball-reference-H1-estimate} implies
	\[
	\| (\xi_r^{-1}\D^{g_r}_{P_{g_r}^\pm}\xi_r)^{-1} - (\D_P^{\mathrm{eucl}})^{-1} \|_{\mathcal{L}(L^2, L^q)} \leq 2 C_E r^2 A_q^{\mathrm{eucl}} B_1^{\mathrm{eucl}}.
	\]
	Accounting for the volume comparison, this implies $A_q(B_1, g_r) = A_q^{\mathrm{eucl}} + O(r^2)$, leading to the asymptotic formula
	\[
	A_q(\Omega_r) = r^{1/(n+1)} \left( A_q^{\mathrm{eucl}} + O(r^2) \right).
	\]
	\end{remark}

\paragraph{\textbf{Application of the fixed-point theorem.}}
Set
\[
p=\frac{2n}{n-1},
\qquad
q=2(p-1).
\]
For each $0<r\leq r_0$, prescribe boundary data through
an extension
$\psi_{0,r}\in H^1(\Sigma\Omega_r)$.
Define its harmonic representative by
\[
\psi_r
=
\psi_{0,r}
-
\D_{P_r}^{-1}(\D\psi_{0,r}).
\]
Then
\[
\D\psi_r=0,
\qquad
P_r\psi_r=P_r\psi_{0,r},
\]
and $\psi_r$ depends only on the prescribed boundary
datum $P_r\psi_{0,r}$.

Assume that, for some constant $M>0$ independent of $r$,
\[
\|\psi_r\|_{L^q(\Omega_r)}
\leq Mr^{n/q}.
\]
Theorem~\ref{thm:geodesic-ball-inverse} gives
\[
\begin{aligned}
	|\lambda|A_q(\Omega_r)\|\psi_r\|_{L^q}^{p-2}
	&\leq
	C|\lambda|M^{p-2}
	r^{1/(n+1)+n(p-2)/q}\\
	&=
	C|\lambda|M^{p-2}r.
\end{aligned}
\]
Therefore, for fixed $\lambda$, $M$, and $\alpha>0$,
the smallness condition
\[
|\lambda|A_q(\Omega_r)\|\psi_r\|_{L^q}^{p-2}
<
\delta_p(\alpha)
\]
holds for all sufficiently small $r$.

For $\psi_r\neq0$,
Corollary~\ref{cor:adapted-radius} consequently gives
an $H^1$ solution of
\[
\D u=\lambda|u|^{p-2}u,
\qquad
P_ru=P_r\psi_r,
\]
unique among solutions satisfying
\[
\|u-\psi_r\|_{L^q(\Omega_r)}
\leq\alpha\|\psi_r\|_{L^q(\Omega_r)}.
\]
The iteration converges geometrically in $L^q$ and,
after the first step, in $H^1$.
For $\psi_r=0$, the corresponding ball consists only
of the zero solution.

\begin{remark}[Boundary data from a fixed smooth section]
	The required estimate on $\psi_r$ holds, in particular,
	when the boundary data are obtained by restricting a
	fixed smooth spinor $\psi_0$ defined near $x_0$.
	Indeed, its harmonic representative on $\Omega_r$ is
	\[
	\psi_r
	=
	\psi_0|_{\Omega_r}
	-
	\D_{P_r}^{-1}
	\bigl((\D\psi_0)|_{\Omega_r}\bigr).
	\]
	Lemma~\ref{lem:harmonic-extension} and
	Theorem~\ref{thm:geodesic-ball-inverse} give
	\[
	\|\psi_r\|_{L^q}
	\leq
	\|\psi_0\|_{L^q(\Omega_r)}
	+
	Cr^{1/(n+1)}
	\|\D\psi_0\|_{L^2(\Omega_r)}.
	\]
	Since $\operatorname{vol}_g(\Omega_r)\leq Cr^n$ and
	\[
	\frac1{n+1}+\frac n2=1+\frac nq,
	\]
	we obtain
	\[
	r^{-n/q}\|\psi_r\|_{L^q}
	\leq
	C\left(
	\|\psi_0\|_{L^\infty(\Omega_r)}
	+
	r\|\D\psi_0\|_{L^\infty(\Omega_r)}
	\right).
	\]
	The right-hand side is uniformly bounded for small $r$.
	Thus, for every fixed $\lambda$, the fixed-point
	theorem applies on all sufficiently small geodesic
	balls with these prescribed boundary data.
\end{remark}
\vspace{0.5cm}
\paragraph{\myuline{\textit{Acknowledgment}}}
\textit{I would like to thank
	my advisor Prof. Dr. Nadine Große for her help and support throughout my studies.}

\printbibliography
\end{document}